\documentclass[leqno,12pt]{article}
\usepackage{hyperref}
\usepackage{latexsym}
 \usepackage{graphicx}
\usepackage{amsmath}
\usepackage{amssymb}
\usepackage{cite}
\usepackage{color}
\usepackage[normalem]{ulem}

\usepackage[margin=1.3in, top=1.3in, bottom=1.3in]{geometry}

\newcommand{\nc}{\newcommand}
\nc{\nt}{\newtheorem}
\nt{thm}{Theorem}[section]
\nt{cor}[thm]{Corollary}
\nt{prop}[thm]{Proposition}
\nt{lem}[thm]{Lemma}
\nt{defn}[thm]{Definition}
\nt{rem}[thm]{Remark}
\nt{exa}[thm]{Example}
\nt{ass}[thm]{Assumption}
\nc{\ip}[2]{\mbox{$\langle #1,#2 \rangle$}}
\nc{\pf}{\noindent{\bf Proof\ \ }}
\nc{\finpf}{\hfill{$\Box$}\linespace}
\nc{\linespace}{\vspace{\baselineskip} \noindent}
\nc{\R}{{\bf R}}
\nc{\E}{{\bf E}}
\nc{\C}{{\bf C}}
\nc{\U}{{\bf U}}
\nc{\M}{{\bf M}}
\nc{\T}{{\bf T}}
\nc{\cA}{{\mathcal A}}
\nc{\cM}{{\mathcal M}}
\nc{\cB}{{\mathcal B}}
\nc{\cS}{{\mathcal S}}
\nc{\cE}{{\mathcal E}}
\nc{\RE}{\mbox{\rm Re}\,}
\nc{\Rn}{{\bf R}^n}
\nc{\Mn}{{\bf M}^n}
\nc{\Cn}{{\bf C}^n}
\nc{\Hn}{{\bf H}^n}

\def\tto{\;{\lower 1pt \hbox{$\rightarrow$}}\kern -12pt
           \hbox{\raise 2.8pt \hbox{$\rightarrow$}}\;}
\newenvironment{myequation}{\setcounter{equation}{\value{thm}}
   \begin{equation}}{\addtocounter{thm}{1}\end{equation}}

\nc{\bmye}{\begin{myequation}}
\nc{\emye}{\end{myequation}}

\begin{document}
\title{
Partial smoothness of the numerical radius at matrices whose fields of values are disks
}
\author{
 A.S. Lewis\thanks{ORIE, Cornell University, Ithaca, NY 14853, U.S.A.
(\texttt{people.orie.cornell.edu/aslewis}).
Research supported in part by National Science Foundation Grant DMS-1613996.}
\and 
M.L. Overton\thanks{Courant Institute of Mathematical Sciences, New York University, New York, NY 10012, U.S.A.
(\texttt{cs.nyu.edu/overton}).
    Research supported in part by National Science Foundation Grant DMS-1620083.}
}
\maketitle

\begin{abstract}
Solutions to optimization problems involving the numerical radius often belong to a special class:  the set of matrices having field of values a disk
centered at the origin. After illustrating this phenomenon with some examples, we illuminate it by studying matrices around which this set of ``disk matrices"
is a manifold with respect to which the numerical radius is partly smooth.  We then apply our results to matrices whose nonzeros consist of a 
single superdiagonal, such as Jordan blocks and the Crabb matrix related to a well-known conjecture of Crouzeix. 
Finally, we consider arbitrary complex three-by-three matrices;  even in this case, the details are surprisingly intricate. 
One of our results is that in this real vector space with dimension 18, the set of disk matrices is a semi-algebraic manifold with dimension~12.

\end{abstract}
%
%

\section{Introduction:  the numerical radius}
Given a matrix $A$ in the space $\Mn$ of $n\times n$ complex matrices, the asymptotic convergence of the discrete-time dynamical system 
\[
x_{k+1} = A x_k \qquad (k=0,1,2,\ldots)
\]
is governed by the spectrum of $A$.  Specifically, the iterates $x_k \in \Cn$ converge to zero for all initial points $x_0$ if and only if $A$ is {\em stable}, meaning its spectrum is contained in the open unit disk.  However, to understand the {\em transient} behavior of the iterates --- how $\|x_k\|$ depends on $k$ --- the spectrum is inadequate.  Instead we can rely on the following well-known inequality:
\bmye \label{power}
\|A^k\|_2 \le 2\big(r(A)\big)^k \qquad (k=0,1,2,\ldots),
\emye
where $\|\cdot\|_2$ denotes the spectral norm (the largest singular value), $r(A)$ denotes the numerical radius of $A$,
\[
      r(A) ~=~ \max \{ |z| : z \in W(A) \},
\]
and the set
\[
       W(A) ~=~ \{u^*Au : u\in\Cn, \|u\|_2 =1 \}
\]
is the field of values (numerical range) of $A$.  For an excellent discussion of this and related issues, see Trefethen and Embree's book \cite[\S14,17]{TreEmb05}.

The numerical radius is a vector space norm, but it does not satisfy the submultiplicative property
required in most definitions of a matrix norm \cite{HorJoh13}. However, Lax and Wendroff, in their famous 1964 paper on difference equations for hyperbolic PDEs \cite{LaxWen64}, 
showed that if $r(A) \leq 1$, then the sequence $(A^k)$ is bounded. Halmos then conjectured the ``power inequality''
\[
      r(A^k) \leq \big(r(A)\big)^k \qquad (k=0,1,2,\ldots).
\]
The first proof, by Berger \cite{Ber65}, was simplified by Pearcy \cite{BerSta67,GolTad82,Lax07,Pea66}.  The inequality (\ref{power}) is an elementary consequence \cite[p.\ 379]{HorJoh13}.
A variety of other properties of the numerical radius are collected in\cite[\S18.4]{HandbookLA}, and a useful bibliography appears in \cite{li}.
Research on the topic remains very active.\footnote{So far this century, at least 150 papers have appeared with titles including the string ``numerical radius" (according to MathSciNet). The concept of numerical radius goes back at least to Lumer (1961) \cite{lumer}, but may well be older.}

An elegant proof of the power inequality relies on the following semidefinite representation of the numerical radius due to \cite{mathias}, based in part on \cite{ando}.  We denote the space of $n$-by-$n$ Hermitian matrices by $\Hn$.

\begin{thm}[Semidefinite representability of numerical radius] \label{mathias} \hfill \mbox{}
The numerical radius of a matrix $A \in \Mn$ is bounded above by $\alpha \in \R$ if and only if there exists a matrix $Z \in \Hn$ such that the Hermitian matrix
\[
\left[
\begin{array}{cc}
\alpha I + Z &  A \\
A^* & \alpha I - Z
\end{array}
\right]
\]
is positive semidefinite. 
\end{thm} 

\noindent
This representation in particular allows us to compute the numerical radius of $A$ by semidefinite programming, choosing variables $\alpha$ and $Z$ to minimize $\alpha$ subject to the semidefinite constraint.  For more direct methods, see \cite{MenOve05,Uhlig2009}.

The numerical radius appears naturally in matrix optimization and best approximation problems: one survey, from an operator-theoretic perspective, appears in \cite{aksoy}.  Our current work is motivated by the following observation: optimization involving the numerical radius quite commonly results in matrices whose fields of values are disks (by which we mean disks in the complex plane centered at zero, including the trivial case $\{0\}$).  We call such matrices {\em disk matrices}, and denote the set of disk matrices by $\cA$.  Surprisingly, an explicit
parametrization of disk matrices was not known until recent work of Crouzeix \cite{Cro16}, discussed below in Section \ref{parametrization}, and based on much older work of Ando \cite{ando}.

Disk matrices are relatively rare:  we focus here primarily on examples when, locally, the set $\cA$  has the structure of a manifold of codimension $2n$ in the $2n^2$-dimensional real vector space $\Mn$.  Our aim here is to illuminate variational properties of the numerical radius (in particular the property of partial smoothness \cite{lewis_active}) that, notwithstanding their rarity, produces disk matrices at optimality.

The paper is organized as follows.  In Section \ref{optimization}, we discuss two interesting examples where disk matrices appear as solutions to optimization problems involving the numerical radius.  The first example involves feedback control, while the second concerns the Crouzeix conjecture in matrix analysis \cite{crouzeix-conjecture,CroPal17}.  Then, in Section \ref{definitions}, we introduce the ideas we use from convex analysis, including the notion of partial smoothness.  Section \ref{examples} summarizes our results for a well-known disk matrix that arises in the study of Crouzeix's conjecture.  In Section \ref{subdifferential} we characterize subgradients of the numerical radius at disk matrices, and then in Sections \ref{unit-disk} and \ref{any-disk} we develop our main results concerning partial smoothness of the numerical radius at disk matrices.  Section \ref{parametrization} discusses our results in the framework of the Crouzeix parametrization of disk matrices, and we end in Section \ref{three-by} with a thorough analysis of the three-by-three case.

\section{Examples of numerical radius optimization} \label{optimization}
For now, we focus on two examples to illustrate our theme.  More examples, arising from the proximal operator for the numerical radius, follow in Section \ref{examples}.

\subsection*{Feedback control}
Our first example generalizes our original dynamical system $x_{k+1} = A x_k$ in the space $\Cn$ by involving an input vector $u_k \in \C^m$ and an output vector $y_k \in \C^p$, related as follows:
\[
            x_{k+1} = Ax_k + B u_k, \qquad y_k= C x_k,
\]
for given matrices $B$ and $C$.  If we set $u_k=Ky_k$ for some matrix $K$ (the controller), we arrive at the system
\[
        x_{k+1} =M(K) x_{k} \qquad (k=0,1,2,\ldots),
\]
defined by the matrix $M(K) = A+BKC$.  The classical static output feedback problem \cite{Won85} is to choose $K$ so that $M(K)$ is stable. The problem of determining whether such a matrix $K$ exists in general seems to be hard;  indeed, when bounds are put on all the entries of $K$, this problem is NP-hard \cite{BloGevLin95}. However, as long as $mp\geq n$ (or $mp > n$ in the
real case), the problem \emph{is} generically solvable via ``pole placement": $K$ can be chosen to place the eigenvalues of
$M(K)$ in any desired location \cite{Wan96, Wil97}.

On the other hand, as we remarked earlier, although stability by itself does not give
information about the transient behavior of the system, such information \emph{is} given by the numerical radius 
$r\big(M(K)\big)$.  By contrast with the question of stability, the problem of minimizing the function 
$r\big(M(\cdot)\big)$ is tractable via
semidefinite programming, thanks to Theorem \ref{mathias}, observing that the map $M$ is affine. 

We now report on some computational experiments minimizing $r\big(M(K)\big)$ which clearly support our
claim that minimizing the numerical radius frequently leads to disk matrices.
Normally $A$, $B$ and $C$ depend on the application of interest, but here we generated them
randomly, choosing them to be real matrices with
i.i.d.\ entries from the normal distribution. 
We set $n=5$ and considered $m$ and $p$ varying from 1 to $n$.
For each choice of $(m,p)$, we generated 1000 random triples $(A,B,C)$ with dimensions
$n\times n$, $n\times m$ and $p\times n$ respectively. 
We then used the CVX package \cite{cvx} (with precision set to ``best") to minimize the numerical radius 
$r\big(M(\cdot)\big)$ over $m$-by-$p$ matrices $K$ for each randomly generated triple, resulting in an
optimal controller $\bar K$. We then used 
the Chebfun package \cite{Chebfun14}, as described in \cite{GreOve17}, to compute the boundary
of the field of values of the optimal matrix $\bar M = M(\bar K)$ and measure how close it is to being a 
circle centered at the origin. To do this we recomputed the numerical radius
\[
r(\bar M)~=~ \max\big\{|z|: z\in \mathrm{bd}\big(W(\bar M)\big)\big\}
\]
as well as the ``inner numerical radius" \cite{CheHig99}
\[
\tilde r(\bar M)~=~ \min\big\{|z|: z\in \mathrm{bd}\big(W(\bar M)\big)\big\},
\]
and computed a measure $r(\bar M) - \tilde r(\bar M)$ of the distortion of $W(\bar M)$ from a disk centered at zero. 
Table 1 shows, for each triple, in how many of the 1000 instances the resulting value $r(\bar M) - \tilde r(\bar M)$
was below a threshold that we set to $10^{-7}$. These quantities are shown as
percentages rounded to the nearest 10\% because repeated runs
indicated that with this choice of rounding, the results are reasonably consistent. For the cases where $0 \%$ is displayed,
actually not even one out of the 1000 runs resulted in a matrix satisfying the criterion $r(\bar M) - \tilde r(\bar M) < 10^{-7}$,
while for the case $m=p=5$ all 1000 optimal matrices satisfied the given criterion; this is simply because in this case, with 25 free
variables, the optimal matrix is zero with probability one. 

\begin{table}
\begin{center}
\begin{tabular}{|c|c|c|c|c|c|}
 \hline & $p=1$ & $p=2$ & $p=3$ & $p=4$ & $p=5$\\ \hline 
$m=1$ & 0 \% & 0 \% & 0 \% & 0 \% & 20 \%\\ \hline 
$m=2$ & 0 \% & 0 \% & 0 \% & 50 \% & 40 \%\\ \hline 
$m=3$ & 0 \% & 0 \% & 70 \% & 70 \% & 70 \%\\ \hline 
$m=4$ & 0 \% & 50 \% & 70 \% & 90 \% & 90 \%\\ \hline 
$m=5$ & 20 \% & 40 \% & 70 \% & 90 \% & 100 \%\\ \hline 
\end{tabular}
\end{center}
\caption{Percentage of cases where $W(M(\bar K))$ is a disk}
\end{table}

More interesting is the fact that these computational results are consistent with our 
theoretical results given later in the paper.  In the cases we study, we
show that the set of disk matrices in the Euclidean space $\Mn$ of complex $n\times n$ matrices
typically has codimension $2n$.  (We do not pursue the real case here, but in the space of $n\times n$ real matrices, the analogous codimension is $n$.)  In the cases in the table where none of the
optimal matrices satisfied the disk matrix criterion, the number of free variables in the optimization problem, $mp$, is
less than $n$, and therefore the probability of the affine space of feasible matrices intersecting the manifold of disk
matrices is zero. In contrast, in all cases when $n \le mp < n^2$, one can expect the affine space of feasible matrices to
intersect the manifold of disk matrices, so that there is at least the possibility that the optimal matrix is a disk matrix ---
and evidently, from the results in the table, the probability of this happening is nonzero, ranging in fact from approximately
20\% to 90\%.

\subsection*{The Crouzeix ratio}
Our second example arises from observations \cite{GreOve17,crouzeix-ratio} that minimizing a certain ratio associated with the Crouzeix conjecture \cite{crouzeix-conjecture,Cro16,CroPal17} seems to result in disk matrices.   In particular, a special case that is known
to be true is that, for matrices 
$X \in \Mn$, the ratio
\bmye \label{cro-ratio}
   \frac{\big(r(X)\big)^{n-1}}{\|X^{n-1}\|_2} 
\emye
is minimized by a bidiagonal matrix whose field of values is a disk. We discuss this matrix, which we call the Crabb matrix, in Section \ref{examples}. 
Note that setting $k=n$ in \eqref{power} implies that the ratio given
above is bounded below by $\frac{1}{2}$; this value is attained by the Crabb matrix.

\section{Subgradients and partial smoothness} \label{definitions}
We begin with some basic convex analysis on Euclidean spaces --- finite-dimensional real inner-product spaces.  We consider several specific Euclidean spaces.  For integers $n \ge 1$, we consider the space 
$\Cn$ of complex column \mbox{$n$-vectors} with inner product 
$\ip{x}{y} = \mbox{Re}(x^*y)$, the space of $n$-by-$n$ complex matrices $\Mn$ with the inner product
$\ip{X}{Y} = \mbox{Re}\,\mbox{trace}(X^*Y)$, and the subspace $\Hn$ of Hermitian matrices (with the same inner product).  Unless otherwise stated, $\|\cdot\|$ always refers to the norm induced by the inner product.  Dimension always refers to the real vector space notion.  A set of vectors $\{x_k : k \in K\}$ is affinely independent in a Euclidean space $\E$ if and only if the vectors $\{(x_k,1) : k \in K\}$ are linearly independent in the product space $\E \times \R$. The notation \mbox{\rm ri} denotes relative interior.

Given a Euclidean space $\E$, and a norm
 $\nu$ on $\E$, a vector $y \in \E$ is, by definition, a subgradient of $\nu$ at a point $x \in \E$ when
\bmye \label{defn}
\ip{y}{w-x} \le \nu(w) - \nu(x) ~~\mbox{for all}~ w \in \E.
\emye
We write the set of subgradients as $\partial\nu(x)$.

Given a nonempty convex set $C \subset \E$, the subspace parallel to $C$ is defined by
\[
\mbox{par}\,C ~=~ \R(C-C).
\]
By contrast, the linear span of $C$ is
\[
\mbox{span}\,C ~=~ \R C - \R C ~\supset~ \mbox{par}\,C.
\]
The two subspaces are identical if and only if $0 \in C$.

\begin{prop} \label{norms}
Consider any norm $\nu$ on a Euclidean space $\E$.  At any point $x \in \E$, we have the following relationships:
\begin{eqnarray*}
\partial \nu (x) &=& \big\{ y \in \partial \nu(0) : \ip{y}{x} = \nu(x) \big\}  \\
\mbox{\rm par}\big(\partial \nu(x)\big) &=& 
x^{\perp} \cap \mbox{\rm span}\big(\partial \nu(x)\big).
\end{eqnarray*}
\end{prop}

\pf
The first equation follows from a standard argument, as follows.  By the definition of subgradients, we have $y \in \partial \nu(0)$ if and only if 
\bmye \label{defn0}
\ip{y}{w} \leq \nu(w) ~~\mbox{for all}~ w \in \E.
\emye
If such a subgradient $y$ also satisfies $\ip{y}{x} = \nu(x)$, then subtracting implies inequality (\ref{defn}).  Conversely, inequality (\ref{defn}) implies $\ip{y}{x} = \nu(x)$, by setting $w$ equal to $2x$ and $0$ in turn, and then adding implies inequality (\ref{defn0}).

Turning to the second equation, any element of the left-hand side can be written in the form
$u = \mu(y-z)$ for a scalar $\mu \in \R$ and subgradients $y,z \in \partial \nu(x)$.  
Clearly $u$ lies in 
$\mbox{\rm span}\big(\partial \nu(x)\big)$, and the first equation guarantees $\ip{u}{x} = 0$, so $u$ also lies in the right-hand side.  

For the converse, since $0 \in \partial \nu(0)$, 
there is nothing to prove when $x=0$, so we can assume 
$x \ne 0$, and hence $\nu(x) \ne 0$.  Any element of the right-hand side can be written in the form  $u = \mu y - \tau z$ for a scalars $\mu,\tau \in \R$ and subgradients $y,z \in \partial \nu(x)$, where
$\ip{\mu y - \tau z}{x} = 0$.  Using the first equation, we deduce $(\mu-\tau)\nu(x) = 0$ so $\mu=\tau$, and the result follows.
\finpf

Following \cite{lewis_active}, we call a norm (or, more generally, an extended-valued convex function) $\nu$ on $\E$ {\em partly smooth} at a point $x \in \E$ relative to a set 
$\cM \subset \E$ containing $x$ if $\cM$ is in fact a \mbox{${\mathcal C}^{(2)}$-smooth} manifold around $x$ whose normal space satisfies \mbox{$N_{\cM}(x) = \mbox{par}\big(\partial \nu(x)\big)$}, and such that the restriction of $\nu$ to $\cM$ is ${\mathcal C}^{(2)}$-smooth and the subdifferential mapping $\partial \nu$ restricted to $\cM$ is continuous.  In particular, any norm is partly smooth at zero relative to $\{0\}$.  This definition captures the smooth/nonsmooth geometry of the norm $\nu$ around the point $x$.  Specifically, the convex function $y \mapsto \nu(x+y)$ restricted to the tangent space
$T_{\cM}(x)$ is differentiable at zero, but is {\em sharp} in nonzero directions  $y \not\in T_{\cM}(x)$:  the sum of the directional derivatives in the directions $y$ and $-y$ is strictly positive.

For example, for any point $x \in \Rn$, defined the index set $J = \{ i : x_j = 0 \}$.  Then the $l_1$-norm $\|\cdot\|_1$ is partly smooth at $x$ relative to the manifold 
\[
\big\{z \in \Rn : z_j = 0~  \mbox{for all}~ j \in J \big\}.
\]
Then the $l_1$-norm is partly smooth at any point in $\cM$ relative to $\cM$.

When a norm is partly smooth relative to a set, we can view that set locally as the set of solutions of a parametrized family of optimization problems.  To illustrate, consider the proximal operator corresponding to a proper closed convex function $\nu$ on $\E$, namely the map 
$\mbox{prox}_{\nu} \colon \E \to \E$ defined by
\[
\mbox{prox}_{\nu}(u) ~=~ 
\mbox{\rm arg}\!\min_{x \in \E} \big\{ \nu(x) + \frac{1}{2}\|x - u\|^2 \big\} \qquad (u \in \E).
\]
Proximal operators are central to many first-order algorithms for convex optimization.  They are often tractable to compute:  in particular, if the function $\nu$ is semidefinite representable, as is the case for the numerical radius (by Theorem \ref{mathias}), then we can compute the corresponding proximal operator by semidefinite programming.

Among several useful consequences of the idea of partial smoothness, particularly intuitive is its use in sensitivity analysis, by which we mean the behavior of solutions to optimization problems under small perturbations to the data.  Under reasonable conditions, the manifold in the definition consists, locally, of all possible optimal solutions to perturbed problems.  

As a typical example, we consider the prox operator associated with a convex function $\nu$.  (See \cite[Cor 4.8]{shanshan}, and \cite{ident}  for more general results along the same lines.)  
Note that $y$ is a subgradient of $\nu$ at a point $x$ if and only if $\mbox{prox}_\nu$ maps $x+y$ to 
$x$.  Partial smoothness reveals what happens under perturbations.

\begin{thm}[Proximal operators and partial smoothness] \label{identification} ~~~~~~\hfill \mbox{}
Suppose a proper closed convex function $\nu$ on a Euclidean space 
$\E$ is partly smooth at a point $x \in \E$ relative to a set $\cM \subset \E$.  Consider any subgradient 
$y \in \mbox{\rm ri}\big(\partial \nu (x)\big)$.  Then 
the proximal operator $\mbox{\rm prox}_\nu$ maps $x+y$ to $x$, and maps 
$x+N$, where $N$ is any small neighborhood of $y$, onto a neighborhood of $x$ in $\cM$.
\end{thm}

\noindent
In particular, with the assumptions of this theorem, with positive probability, the proximal operator 
$\mbox{\rm prox}_{\nu}$ maps vectors chosen randomly from a probability distribution with support $\E$ into the set $\cM$.  Thus optimization problems (and proximal operators) commonly produce solutions in the manifold $\cM$, even though its dimension may be low.  The results we prove below about partial smoothness of the numerical radius relative to the set of disk matrices suggest that analogous behavior may be at play in the numerical examples we presented in Section \ref{optimization}. 

\section{The Crabb matrix}  \label{examples}
In this work, we prove partial smoothness of the numerical radius around some interesting matrices.  Before proceeding to our general development, we illustrate our results on some simple special cases. Here and throughout the paper, we denote the unit circle in the complex plane by $\T$.

\subsubsection*{The two-by-two case}
In ${\mathbf M}^2$, disk matrices coincide with nilpotent matrices, as follows immediately from \cite{Cro16}.  In particular, consider
the Jordan block
\[
J = 
\left[
\begin{array}{cc}
0 & 2 \\
0 & 0 
\end{array}
\right]
\]
scaled so that its field of values is the unit disk. 
As shown in \cite[p.~231 and p.~238]{crouzeix-ratio}, at $J$, the numerical radius $r$ has subdifferential
\[
\partial r(J) ~=~ \frac{1}{2} \mbox{conv} \Big\{ 
\left[
\begin{array}{cc}
w & 1 \\
w^2 & w 
\end{array}
\right]
: w \in \T \Big\},
\]
and furthermore, integrating over $w \in \T$ shows
$\frac{1}{4} J \in \mbox{ri} \big( \partial r(J) \big)$ (see \cite[p.~242]{crouzeix-ratio}).
As we show later, $r$ is partly smooth at $J$ relative to the set of disk  matrices $\cA$, or equivalently in the case $n=2$, the set of nilpotent matrices.  Hence, by Theorem \ref{identification}, the proximal operator for the numerical radius, $\mbox{prox}_r$, maps all matrices near $\frac{5}{4} J$ to nilpotent matrices near $J$.  Computational tests, using semidefinite programming, confirm this property.

\subsubsection*{The three-by-three case}
The analogous behavior in ${\mathbf M}^3$ is also simple to check.  Now we consider the Jordan block
\[
J = 
\left[
\begin{array}{ccc}
0 & \sqrt{2} & 0\\
0 & 0 & \sqrt{2}\\
0 & 0 & 0 
\end{array}
\right],
\]
again scaled so that its field of values is the unit disk.
Following the same argument as in \cite[p.~231]{crouzeix-ratio},
it is easy to see that at $J$, 
the numerical radius $r$ has subdifferential
\[
\partial r(J) ~=~\frac{1}{4} \mbox{conv} \Big\{ 
\left[
\begin{array}{ccc}
w 			& \sqrt{2} 		&  w^{-1}	\\
\sqrt{2}w^2	& 2w			& \sqrt{2}			\\
 w^3	& \sqrt{2}w^2	& w  
\end{array}
\right]
: w \in \T \Big\}.
\]
Integrating again, we deduce $\frac{1}{4} J \in \mbox{ri} \big( \partial r(J) \big)$.
As we show later, $r$ is partly smooth at $J$ relative to $\cA$, so the proximal operator 
$\mbox{prox}_r$ maps all matrices near $\frac{5}{4} J$ to disk matrices near $J$.  
Computational tests again confirm this property.

As we shall see (Theorem \ref{3disk}), a necessary condition for $X \in\M^3$ to be a disk matrix is:  
\begin{quote}
$X$ is unitarily similar to a matrix of the form 
\[
\left[
\begin{array}{ccc}
0 & a & b \\
0 & c & d \\
0 & 0 & 0
\end{array}
\right]
\]
for numbers $a,b,c,d \in \C$ satisfying the condition $c(|a|^2 + |d|^2) = -adb^*$.
\end{quote}
In fact this condition along with the condition that $a$ and $d$ are both nonzero
is also sufficient for $X \in \cA$ (by the same result).  For $X$ near $J$, these necessary and sufficient conditions are equivalent, and furthermore easy to check computationally by reduction to Schur form.  We show (Theorem \ref{three-final}) that the numerical radius is partly smooth at all such matrices, relative to $\cA$.

By contrast, we argue at the end of this work (Proposition \ref{failure}) that the numerical radius $r$ is {\em not} partly smooth relative to the set of disk matrices $\cA$ (or indeed any subset) at the disk matrix
\[
E_0 ~=~ 
\left[
\begin{array}{ccc}
0 & 0 & 1 \\
0 & 0 & 0 \\  
0 & 0 & 0
\end{array}
\right].
\]
Our later results imply
$\frac{1}{2} E_0 ~\in~ \mbox{ri}\big( \partial r(E_0) \big)$,
so certainly we have $\mbox{prox}_r\big(\frac{3}{2}E_0\big) = E_0$.
However, the lack of partial smoothness means that there is no guarantee that the proximal operator
$\mbox{prox}_r$ maps matrices near $\frac{3}{2}E_0$ to disk matrices.
%

\subsubsection*{The general case}
Extending the three-by-three case to $\Mn$ for $n > 3$ involves not a Jordan block but the Crabb matrix
\cite{Cra71,GusRao97}\footnote{This matrix was called the ``Choi-Crouzeix" matrix in \cite{crouzeix-ratio,GreOve17}.
Many thanks to A.~Salemi for providing the reference \cite{Cra71} and pointing out that as a consequence of Crabb's results, 
the conjectures on \cite[p.~239]{GreOve17} and \cite[p.~240]{GreOve17} are true. A third conjecture on \cite[p.~242]{GreOve17} remains open.}
\[
     K ~=~ \left [ \begin{array}{ccccccc}
                 0 & \sqrt{2} &        &       &       &      &           \\
                         & \cdot    & 1      &       &       &      &           \\
                         &          & \cdot  & \cdot &       &      &           \\
                         &          &        & \cdot & \cdot &      &           \\
                         &          &        &       & \cdot & 1    &           \\
                         &          &        &       &       & \cdot & \sqrt{2} \\
                         &          &        &       &       &       &  0 \end{array} \right ]  .
\]
This matrix has field of values the unit disk, and again following the reasoning in \cite[p.~231]{crouzeix-ratio}, the matrix
\[
H ~=~ \frac{1}{n-1}
     \left [ \begin{array}{ccccccc}
                 0 & \frac{1}{\sqrt{2}} &        &       &       &      &           \\
                         & \cdot    & 1      &       &       &      &           \\
                         &          & \cdot  & \cdot &       &      &           \\
                         &          &        & \cdot & \cdot &      &           \\
                         &          &        &       & \cdot & 1    &           \\
                         &          &        &       &       & \cdot & \frac{1}{\sqrt{2}} \\
                         &          &        &       &       &       &  0 \end{array} \right ]
\]
is a subgradient of the numerical radius $r$ at $K$, in fact belonging to 
$\mbox{ri}\big(\partial r(K)\big)$.  As we show later, $r$ is partly smooth at $K$ relative to the disk matrices $\cA$, so the proximal operator $\mbox{prox}_r$ maps all matrices near $K+H$ to disk matrices near $K$.  
Again, we can verify this behavior computationally.

As we shall see, a necessary condition for a matrix $X$ to be a disk matrix (by Proposition \ref{eigenvalue}) is that the maximal eigenvalues $\lambda_{\max}(w^*X + wX^*)$, as $w$ ranges over any subset of $\T$ of cardinality $2n+1$, are all equal.  This necessary condition is easy to check computationally, and in fact it is also sufficient for $X$ near $K$ (by Proposition \ref{check}).

\section{The subdifferential of the numerical radius at disk matrices} \label{subdifferential}
Clearly both the field of values $W(X)$ and the numerical radius $r(X)$ are invariant under unitary similarity transformations on the matrix $X$.  Since we can write
\bmye \label{numerical}
r(X) ~=~ \max \big\{ \mbox{Re}(w^*u^* X u) : w \in \C,~ |w|=1,~ u \in \Cn, \|u\|=1 \big\},
\emye
it follows immediately that the function $r$ is a norm.  

Our analysis depends heavily on the ``Hermitian part'' linear map $\Phi \colon \Mn \to \Hn$, defined by
\[
\Phi(X) = \frac{1}{2}(X + X^*) \qquad (X \in \Mn).
\]
Notice that this map is the adjoint of the natural embedding:  $\Phi^* \colon \Hn \to \Mn$ is given by
\bmye \label{adjoint}
\Phi^*(Z) = Z \qquad (Z \in \Hn),
\emye
since
\begin{eqnarray*}
\lefteqn{
\big<X , \Phi^*(Z) \big> ~=~ \big< \Phi(X) , Z \big> ~=~ \Big< \frac{1}{2}(X + X^*) , Z \Big>
~=~ \frac{1}{2}\big(\ip{X}{Z} + \ip{X^*}{Z}\big) 
} \\
& & 
=~ \frac{1}{2}\big(\mbox{Re}\, \mbox{tr}(X^*Z) + \mbox{Re}\, \mbox{tr}(XZ)\big)
~=~ \frac{1}{2}\big(\mbox{Re}\, \mbox{tr}(X^*Z) + \mbox{Re}\, \mbox{tr}(ZX^*)\big)  \\
& & 
=~ \frac{1}{2}\big(\mbox{Re}\, \mbox{tr}(X^*Z) + \mbox{Re}\, \mbox{tr}(X^*Z)\big)
~=~ \ip{X}{Z}.
\end{eqnarray*}

We will be interested in the largest eigenvalue $\lambda_{\max}$ of matrices of the form $\Phi(w^*X)$ for scalars $w$ on the unit circle 
$\T \subset \C$.  We use the following result (dating back to \cite{kippenhahn}) extensively:  since its proof is simple, we include it.

\begin{prop}[Disk matrix characterization] \label{eigenvalue}
For any matrix $X \in \Mn$,
\[
\max \big\{ \ip{w}{z} : z \in W(X) \big\} ~=~
\lambda_{\max}\big(\Phi(w^*X)\big) ~~\mbox{for all}~ w \in \T.
\]
Hence $X$ has numerical radius
\[
r(X) = \max_{w \in \T}  \lambda_{\max} \big( \Phi(w^*X) \big),
\]
and has field of values a disk of radius $\rho \ge 0$ if and only if
$\lambda_{\max}\big(\Phi(w^*X)\big) = \rho$
for all $w \in \T$. 
\end{prop}

\pf
Consider the first equation.  The field of values is nonempty and compact, so the maximum on the left-hand side is indeed attained.  Furthermore we have
\begin{eqnarray*}
\lefteqn{
\max \big\{ \ip{w}{z} : z \in W(X) \big\}
} \\
 &=& 
\max \big\{ \mbox{Re}(w^*u^* X u) : u \in \Cn, \|u\|=1 \big\} \\
 &=& 
\max \big\{ \frac{1}{2}(w^*u^* X u + u^* X^* uw) : u \in \Cn, \|u\|=1 \big\} \\
 &=&
\max \Big\{ u^*\big( \frac{1}{2}(w^*X + wX^*)\big)u : u \in \Cn, \|u\|=1 \Big\} \\
 &=&
\lambda_{\max}\big(\Phi(w^*X)\big),
\end{eqnarray*}
so the first equation follows.  The second equation is an immediate consequence.  One direction in the final equivalence is also immediate, and the converse follows since the field of values $W(X)$ is convex by the Toeplitz-Hausdorff Theorem.
\finpf

\begin{lem} \label{gradient}
Given a matrix $\hat X \in \Mn$ and a scalar $w \in \T$, suppose that the matrix 
$\Phi(w^* \hat X)$ has simple largest eigenvalue, with a corresponding unit eigenvector $g$.  Then
the function $X \mapsto \lambda_{\max}\big(\Phi(w^*X)\big)$ is analytic on $\Mn$ around $\hat X$, and its gradient there is $wgg^*$.  Furthermore, the following condition holds:
\[
\ip{wgg^*}{\hat X} ~=~ \lambda_{\max}\big(\Phi(w^*\hat X)\big).
\]
\end{lem}

\pf
Since the largest eigenvalue of the matrix $\Phi(w^* \hat X)$ is simple, it is well known  \cite{tsing} that the function $\lambda_{\max} \colon \Hn \to \R$ is locally analytic there, with gradient $gg^*$.  Therefore, for small matrices $Z \in \Hn$ we have
\begin{eqnarray*}
\lefteqn{
\lambda_{\max}\Big(\Phi\big(w^*(\hat X+Z)\big)\Big) - \lambda_{\max}\big(\Phi(w^*\hat X)\big)
} \\
& &  =
\lambda_{\max}\Big(\Phi\big(w^*\hat X) + \Phi(w^*Z)\big)\Big) - \lambda_{\max}\big(\Phi(w^*\hat X)\big)  \\
& & = \big<g g^* , \Phi(w^*Z) \big> + o(\|Z\|).
\end{eqnarray*}
Using the adjoint property (\ref{adjoint}), notice
\[
\big<g g^* , \Phi(w^*Z) \big> ~=~ \big<\Phi^*(g g^*) , w^*Z \big> ~=~
\mbox{Re}\, \mbox{tr} \big( g g^* w^*Z \big)
~=~ \big< wg g^* , Z \big>.
\]
The expression for the gradient follows.  We then note
\begin{eqnarray*}
\lefteqn{ 
\big< wg g^* , \hat X \big> ~=~
\RE \mbox{tr} \big(w^*g g^* \hat X \big) ~=~ \RE\! \big(g^*w^*\hat X g \big) 
} \\
& & \hspace{2cm} =~ g^* \big(\Phi(w^*\hat X)\big) g ~=~ g^* \lambda_{\max}\big(\Phi(w^*{\hat X})\big) g ~=~ 
\lambda_{\max}\big(\Phi(w^*{\hat X})\big),
\end{eqnarray*}
using the definition of $g$.
\finpf

We next define a set-valued mapping $\Delta \colon \Mn \tto \Mn$ by 
\bmye  \label{G}
\Delta(X) ~=~  
\big\{ wgg^* : w \in \T,~ \|g\|=1,~ \Phi(w^*X) g =  \lambda_{\max}\big( \Phi(w^*X) \big) g \big\},
\emye
and a set of matrices
\bmye
\Omega ~=~
\big\{ X \in \Mn : \lambda_{\max}\big( \Phi(w^*X) \big)~ \mbox{\rm is simple for all}~ w \in \T \big\}
\label{Omega}  
\emye

\begin{prop} \label{continuity}
On the open set $\Omega$, the set-valued mappings $X \mapsto \Delta(X)$ and 
$X \mapsto \mbox{\rm conv}\big( \Delta(x) \big)$ are both continuous.
\end{prop}

\pf
Suppose the set $\Omega$ is not open, so there exists a matrix $X \in \Omega$, and a sequence $X_r \to X$ with $X_r \not\in \Omega$ for all $r=1,2,3,\ldots$.  By definition, for each $r$ there exists a scalar $w_r \in \T$ such that the second largest eigenvalue $\lambda_2 \big( \Phi(w_r^*X_r) \big)$ equals the largest $\lambda_{\max}\big( \Phi(w_r^*X_r) \big)$.  After taking a subsequence, we can suppose $w_r$ converges to some scalar $w \in T$, but then, since the map $\Phi$ and both the eigenvalues 
$\lambda_{\max}$ and 
$\lambda_2$ are continuous, we deduce the contradiction
\[ 
\lambda_2 \big( \Phi(w^* X) \big) = \lambda_{\max}\big( \Phi(w^* X) \big).
\]
Hence $\Omega$ is indeed open.  

To prove continuity, we use Lemma \ref{gradient}, which shows that, for any scalar $w \in \T$, the function $F \colon \Mn \to \R$ defined by $F(X) = \lambda_{\max}\big( \Phi(w^* X) \big)$ is analytic throughout the set $\Omega$.  Furthermore, the gradient of $F$ at any matrix $X \in \Omega$ is given by 
$\nabla F(X) = wgg^*$ for any unit eigenvector $g$ corresponding to the eigenvalue
$\lambda_{\max}\big( \Phi(w^* X) \big)$.  Now consider a sequence $X_r \to X$ in $\Omega$.  We prove $\Delta(X_r) \to \Delta(X)$ in two steps.  

First, we observe that any limit point $Y$ of a sequence $Y_r \in \Delta(X_r)$ must lie in $\Delta(X)$.  To see this, after taking a subsequence we can suppose 
$Y_r \to Y$, and $Y_r = w_r g_r g_r^*$ for some scalar $w_r \in \T$ and unit vector $g_r$ satisfying 
\[
\Phi(w_r^* X_r) g_r ~=~ \lambda_{\max}\big( \Phi(w_r^* X_r) \big) g_r ~~~\mbox{for all}~r.
\]
Taking a further subsequence, we can suppose $w_r$ approaches some scalar $w \in \T$ and $g_r$ approaches some unit vector $g$.  By continuity we deduce
\[
Y = wgg^* ~~~\mbox{and}~~~ \Phi(w^* X) g ~=~ \lambda_{\max}\big( \Phi(w^* X) \big) g,
\]
so we have shown $Y \in \Delta(X)$, as desired.

Secondly, given any matrix $Y \in \Delta(X)$, we prove the existence of a sequence $Y_r \in \Delta(X_r)$ converging to $Y$.  To do so, first write $Y = wgg^*$ for some unit eigenvector $g$ corresponding to the eigenvalue
$\lambda_{\max}\big( \Phi(w^* X) \big)$, so $Y = \nabla F(X)$.  If we define $Y_r = \nabla F(X_r)$ for each $r$, then analyticity implies $Y_r \to Y$, and applying Lemma \ref{gradient} at each $X_r$ shows $Y_r \in \Delta(X_r)$, as desired.

Continuity of the map $\Phi$ and the eigenvalue $\lambda_{\max}$ shows that, for any matrix 
$X \in \Omega$, the image $\Delta(X)$ is always closed, and clearly every element has norm one, so the mapping $\Delta$ is also locally (in fact globally) bounded.  Continuity of the mapping $X \mapsto \mbox{\rm conv}\big( \Delta(x) \big)$ now follows in a standard fashion \cite[Prop.\ 4.30]{VA}.
\finpf 

\begin{thm}[Disk matrix subgradients] \label{subdifferential-continuity}
At any matrix $X$ in the set $\cA$ of matrices whose fields of values are disks centered at zero, the subdifferential of the numerical radius $r$ satisfies
\begin{eqnarray*}
\partial r(X) 
&=& 
\mbox{\rm conv}\big(\Delta(X)\big) \\
\mbox{\rm par}\big(\partial r(X) \big) 
&=&
X^{\perp} \cap \mbox{\rm span}\big(\Delta(X)\big),
\end{eqnarray*}
where the set-valued mapping $\Delta$ is defined by equation (\ref{G}).  On the intersection of $\cA$ with the set $\Omega$ defined by equation (\ref{Omega}), the numerical radius is analytic and its subdifferential is continuous.
\end{thm}

\pf
We can write equation (\ref{numerical}) as 
\begin{eqnarray*}
r(X) 
&=&  
\max_{w \in \T} \max_{\|g\|=1} \ip{wgg^*}{X},
\end{eqnarray*}
for all matrices $X \in \Mn$, so 
\[
\partial r(X) ~=~ \mbox{conv} \big\{wgg^* : w \in \T,~ \|g\|=1,~ \ip{wgg^*}{X} = r(X) \big\}
\]
(see \cite[10.31]{VA}).
For any fixed scalar $w \in \T$, the inner maximum in the expression for the numerical radius equals the eigenvalue 
$\lambda_{\max}\big( \Phi(w^* X) \big)$, and the attaining set  consists of all corresponding unit eigenvectors $g$.  For any matrix $X \in \cA$, the outer maximum is attained on the whole set $\T$, and hence 
$\partial r(X) = \mbox{conv}\big(\Delta(X)\big)$.  The second equation follows from Proposition \ref{norms}.
Continuity follows by Proposition \ref{continuity}, and since for all $X \in \cA$ we have
$r(X) = \lambda_{\max}\big( \Phi(X) \big)$, analyticity follows from Lemma \ref{gradient}.
\finpf

\section{Matrices with field of values the unit disk} \label{unit-disk}
Let $\cB$ denote the set of matrices with field of values the \emph{unit} disk; hence $\cA=\R_+ \cB$.
In this section we describe conditions under which $\cB$ is locally a manifold.  We focus on the subset $\cB \cap \Omega$, consisting of those matrices $X \in\cB$ for which the largest eigenvalue of the matrix 
$\Phi(w^* X)$ equals one and is simple for all scalars $w \in \T$.  We begin with a tool for characterizing when a matrix is
in $\cB$, or equivalently, in $\cB \cap \Omega$, given that it is close to a given matrix in 
$\cB \cap \Omega$.

\begin{prop} \label{equivalent}
Consider any matrix $\hat X$ in the set $\cB  \cap \Omega \subset \Mn$.  Then
for any matrix $X$ near $\hat X$, and any set $\hat\T \subset \T$ of cardinality at least $2n+1$, the following properties are equivalent.
\begin{enumerate}
\item[{\rm (i)}]
$X \in \cB$.
\item[{\rm (ii)}]
$\lambda_{\max} \big( \Phi(w^*X) \big) = 1$ for all $w \in \hat\T$.
\item[{\rm (iii)}]
$\Phi(w^* X) - I$ is singular for all $w \in \hat\T$. 
\item[{\rm (iv)}]
$\Phi(w^* X) - I$ is singular for all $w \in \T$. 
\item[{\rm (v)}]
$\lambda_{\max} \big( \Phi(w^*X) \big) = 1$ for all $w \in \T$.
\end{enumerate}
\end{prop}

\pf
The implications (v) $\Rightarrow$ (i) $\Rightarrow$ (ii) follow from Proposition \ref{eigenvalue}.  The implication (ii) $\Rightarrow$ (iii) is immediate.  

Notice the identity
\[
w^n \det \big( \Phi(w^* X) - I \big) ~=~ \det \Big( \frac{1}{2}(X + w^2 X^*) - wI \Big)
~~~ \mbox{for all $w \in \T$}.
\]
The right-hand side is a polynomial in $w$ of degree at most $2n$, so it vanishes everywhere on the set 
$\hat\T$ if and only if it vanishes everywhere on $\T$.  It follows that properties (iii) and (iv) are equivalent.

Suppose the implication (iv) $\Rightarrow$ (v) fails.  Then
there exists a sequence of matrices $X_r$ converging to $\hat X$, and scalars $w_r \in \T$, such that $1$ is an eigenvalue of the matrix $\Phi(w_r^*X_r)$ and yet is not its largest eigenvalue.  Consequently, its second largest eigenvalue $\lambda_2$ satisfies $\lambda_2 \big( \Phi(w_r^*X_r) \big) \ge 1$ for all $r$.  After taking a subsequence, we can suppose $w_r \to w \in \T$, and we deduce by continuity 
$\lambda_2 \big( \Phi(w^*\hat X) \big) \ge 1$.  But this contradicts our assumption that the matrix
$\Phi(w^*\hat X)$ has a simple largest eigenvalue equal to $1$.
\finpf

We now define functions $f \colon \T \to \Cn$ and $F \colon \T \to \Mn$ by
\[
f(w) = 
\left[
\begin{array}{c}
w \\ w^2 \\ w^3 \\ \vdots \\ w^n
\end{array}
\right]
\]
and
\[
F(w) ~=~ wf(w)f(w)^*
~=~
\left[
\begin{array}{cccccc}
w		& 1 		& w^{-1} 	& w^{-2} 	& \cdots 	& w^{2-n}	\\
w^2		& w 		& 1			& w^{-1} 	& \cdots 	& w^{3-n}	\\
w^3		& w^2 		& w			& 1		 	& \cdots 	& w^{4-n}	\\
w^4		& w^3 		& w^2		& w		 	& \cdots 	& w^{5-n}	\\
\vdots	& \vdots	& \vdots	& \vdots 	& \ddots 	& \ddots	\\
w^n		& w^{n-1} 	& w^{n-2} 	& w^{n-3} 	& \cdots 	& w
\end{array}
\right]
\]

\begin{prop} \label{li}
For any set $\tilde\T \subset \T$ of cardinality at most $2n+1$, the set of vectors 
$\{f(w) : w \in \tilde\T \}$
is affinely independent in the Euclidean space $\Cn$, and the set of matrices
$\{F(w) : w \in \tilde\T \}$ 
is linearly independent in the Euclidean space $\Mn$.
\end{prop}

\pf
We begin with the first claim. Suppose that it fails for some set $\tilde\T$ of distinct scalars $w_1,w_2,w_3,\ldots,w_m$ where 
$m \le 2n+1$.  Without loss of generality, we can assume $m = 2n+1$:  otherwise, simply augment the set 
$\tilde\T$.

Then the linear map $\Psi \colon \R^{2n+1} \to \R \times \Cn$ defined by
\[
\Psi(y) ~=~ \left( \sum_j y_j , \Big(\sum_{j=1}^{2n+1} y_j w_j^k\Big)_{k=1}^n \right)   
~~~(y \in \R^{2n+1})
\]
is not one-to-one.  Since $\dim \R^{2n+1} = 2n+1 = \dim(\R \times \Cn)$, the map $\Psi$ is not onto either, and hence its adjoint $\Psi^* \colon \R \times \Cn \to \R^{2n+1}$ is not one-to-one.  In other words, there exists a scalar $\alpha \in \R$ and $x \in \Cn$, not both zero, such that $\Psi^*(\alpha,x) = 0$.  By definition, for all vectors $y \in \R^{2n+1}$ we have
\begin{eqnarray*}
\lefteqn{
\big< \Psi^*(\alpha,x) , y \big> ~=~ \big< (\alpha,x) , \Psi(y) \big>
} \\
& & 
=~ \alpha \sum_{j=1}^{2n+1} y_j + \mbox{Re}\Big( \sum_{k=1}^n x_k^* \sum_{j=1}^{2n+1} y_j w_j^k \Big)
~=~ \sum_{j=1}^{2n+1} y_j \Big(\alpha + \mbox{Re}\sum_{k=1}^n x_k^* w_j^k \Big)
\end{eqnarray*}
so for $j=1,2,3,\ldots,2n+1$ we have $0 = \big( \Psi^*(\alpha,x) \big)_j =  
\alpha + \mbox{Re}\sum_k x_k^* w_j^k$.
However, the nonzero real trigonometric polynomial 
$\theta \mapsto \alpha + \mbox{Re}\sum_k x_k^* e^{ik\theta}$ can have at most $2n$ zeroes \cite[p.\ 150]{powell-book}, which is a contradiction.

This completes the proof of the first claim.  The second claim follows immediately, since, for all 
$w \in \T$, the first column of the matrix $F(w)$ is the vector $f(w)$, and the $(1,2)$ entry equals $1$.
\finpf

\begin{cor} \label{dimension}
In the Euclidean space $\C^{n+1}$, the convex set
\[
\mbox{\rm conv} 
\Big\{ 
\left[
\begin{array}{c}
1 \\ f(w) 
\end{array}
\right]
: w \in \T 
\Big\}
\]
has dimension $2n$.
\end{cor}

\pf
The given set, which we denote by $S$, contains the set  
\[
\Big\{ 
\left[
\begin{array}{c}
1 \\ f(w) 
\end{array}
\right]
: w \in \hat\T 
\Big\}
\]
for any set $\hat\T \subset \T$.  By choosing $\hat\T$ of cardinality exactly $2n+1$, we see by Proposition \ref{li} that $S$ contains $2n+1$ linearly independent vectors, and hence has dimension at least $2n$.  On the other hand, $y_1 = 1$ for all vectors $y \in S$, so $S$ can have dimension no more than $2n$.
\finpf

To prove that the set $\cB$ of matrices with field of values the unit disk is a manifold locally, we focus on matrices $X \in \cB  \cap \Omega$, and make one further assumption:  that the matrix $\Phi(w^*X)$ parametrized by $w \in \T$ has an (unnormalized) eigenvector corresponding to the simple eigenvalue
$\lambda_{\max} \big( \Phi(w^*X) \big) = 1$ whose components are linearly independent polynomials in $w$ of degree less than $n$.  At first sight, this assumption seems artificial, but in the sections that follow we provide motivation and examples.

\begin{thm} \label{normal}
Consider any matrix $X \in \Mn$.  Suppose there exists an invertible matrix $G \in \Mn$ such that the largest eigenvalue of the matrix $\Phi(w^* X)$ is simple and equal to one for all scalars $w \in \T$, with a corresponding eigenvector $Gf(w)$.
Then $X\in\cB$ and the set $\cB$ of matrices with field of values the unit disk is an analytic manifold around $X$ of codimension $2n+1$.  Its tangent and normal spaces at $X$ are
\begin{eqnarray*}
T_{\cB}(X) &=& \Big\{ Y \in \Mn : \ip{GF(w)G^*}{Y} = 0~ \mbox{for all}~ w \in \T \Big\} \\
N_{\cB}(X) &=& G \Big(\mbox{\rm span}\big\{ F(w) : w \in \T \big\}\Big) G^*.
\end{eqnarray*}
\end{thm}

\pf
That $X\in\cB$ is immediate from Proposition \ref{equivalent}.
For each scalar \mbox{$w \in \T$}, denote the normalized eigenvector $\|Gf(w)\|^{-1}Gf(w)$ by $g(w)$.  
Choose any any set $\hat\T \subset \T$ of cardinality exactly $2n+1$.  
By Proposition~\ref{equivalent}, the set $\cB$ is described locally by the analytic equations in its property (ii).  We see from Lemma~\ref{gradient} that the gradients of the left-hand sides of these equations are linearly independent, since for any $y \in \R^{2n+1}$, the relationship
\[
\sum_j y_j w_j g(w_j) g(w_j)^* = 0
\]
implies 
\begin{eqnarray*}
0 
&=& 
\sum_j y_j w_j \|Gf(w_j)\|^{-2} G f(w_j) f(w_j)^*G^* \\
&=& 
G \Big(\sum_j y_j \|Gf(w_j)\|^{-2} \big(w_j f(w_j) f(w_j)^*\big) \Big) G^* \\
&=& 
G \Big(\sum_j y_j \|Gf(w_j)\|^{-2} F(w_j) \Big) G^*,
\end{eqnarray*}
and hence
\[
\sum_j \big(y_j \|Gf(w_j)\|^{-2}\big) F(w_j) = 0,
\]
and then Proposition \ref{li} shows  $y_j \|Gf(w_j)\|^{-2} = 0$ for all $j$ and so
$y=0$.  The manifold property follows and
\[
T_{\cB}(X) ~=~ \Big\{ Y \in \Mn : \ip{wg(w)g(w)^*}{Y} = 0~ \mbox{for all}~ w \in \hat\T \Big\}.
\]
Since $\hat\T$ was arbitrary, the formula for the tangent space follows, and that for the normal space is an immediate consequence.
\finpf

\section{Disk matrices} \label{any-disk}
Our next step is to lift the main result in the previous section from the set $\cB$ to the set of disk matrices $\cA$.  We first note a computational technique for recognizing 
when a matrix is in $\cA$, given that it is close to a given matrix in $\cA \cap \Omega$.

\begin{prop}[Recognizing disk matrices] \label{check}
Consider any fixed matrix $\hat X$ in the set $\cA \cap \Omega \subset \Mn$ and any set $\hat\T \subset \T$ of cardinality at least $2n+1$.  A necessary condition for a matrix $X \in \Mn$ to be a disk matrix is that the eigenvalues $\big\{ \lambda_{\max} \big( \Phi(w^*X) \big) : w \in \hat\T \big\}$ should all be equal.  
For $X$ near $\hat X$, this condition is also sufficient.
\end{prop}

\pf
The case $n=1$ is an easy exercise, relying on the fact that a scalar $x \in \C$ is zero if and only if  $\RE\!(w^*x)$ takes identical values for at least three distinct scalars $w \in \T$.  We can therefore assume $n>1$.

The necessary condition is clear from Proposition \ref{eigenvalue}.  Suppose the converse fails, so there exists a sequence of matrices $(X_s)$ outside $\cA$ converging to $\hat X$ and a real sequence 
$(\alpha_s)$ such that
\[
\lambda_{\max} \big( \Phi(w^*X_s) \big) = \alpha_s ~~\mbox{for all}~ w \in \hat\T,~ \mbox{and all}~s.
\]
Clearly we have $\alpha_s \to r(\hat X) > 0$.  In particular, we can assume $\alpha_s > 0$ for all 
$s$.  But now we have
\[
\frac{1}{\alpha_s} X_s ~\to~  \frac{1}{r(\hat X)} \hat X ~\in~ \cB \cap \Omega,
\]
so we can apply Proposition \ref{equivalent} to deduce the contradiction $\frac{1}{\alpha_s} X_s \in \cB$ for all large $s$.
\finpf

We next establish a simple tool.

\begin{lem} \label{cone}
Consider any matrix $\hat X \in \cB$.
For all real $\delta > 0$, the set $\cA$ of disk matrices coincides with the set
\[
\cA_{\delta} ~=~ \big\{ \rho X : X \in \cB,~ \rho \in \R,~ 
\|X - \hat X\| < \delta,~ |\rho-1| < \delta \big\}
\]
around $\hat X$.
\end{lem}

\pf
Since the field of values satisfies $W(\rho X) = \rho W(X)$ for all real $\rho > 0$, we clearly always have
$\cA_{\delta} \subset \cA$.  Conversely, consider any sequence of matrices $X_k \to \hat X$ in $\cA$.  The numerical radius $r$ is continuous, so $r(X_k) \to r(\hat X) = 1$.  Now note 
$\frac{1}{r(X_k)}X_k \in \cB$, so 
$X_k \in \cA_{\delta}$ eventually.  The result follows.
\finpf

\begin{prop} \label{tool}
Consider a Euclidean space $\E$, and a set $\cS \subset \E$ that is a \mbox{$k$-dim}\-ensional analytic manifold around a point $\bar s \in \cS$.  Suppose that $\bar s$ does not lie in the tangent space $T_{\cS}(\bar s)$.  Then for all small $\delta > 0$, the set
\[
\big\{ \rho s : s \in \cS,~ \rho \in \R,~ \|s - \bar s\| < \delta,~ |\rho-1| < \delta \big\}
\]
is a $(k+1)$-dimensional manifold around $\bar s$, and its tangent space there is 
$T_{\cS}(\bar s) + \R\bar s$.
\end{prop}

\pf
Denote the open unit ball in the space $\E$ by $B$.  By choosing suitable local coordinates for the manifold $\cS$ we can suppose there is an analytic map $P \colon B \to \E$ such that $P(0)=\bar s$ and with derivative $DP$ one-to-one throughout $B$ and such that the image $P(B)$ is an open neighborhood of 
$\bar s$ in $\cS$.  In that case the tangent space $T_{\cS}(\bar s)$ is just the range of $DP(0)$.

Now define an analytic map $Q \colon B \times \R \to \E$ by $Q(y,\rho) = \rho P(y)$.  Notice
\[
0 = DQ(0,1)(y,\rho) = DP(0)y + \rho \bar s 
\]
implies $(y,\rho) = (0,0)$, since by assumption the vector $\bar s$ does not lie in the range of $DP(0)$ and the map $DP(0)$ is one-to-one.  Therefore the derivative map $DQ$ is also one-to-one at the vector $(0,1)$, so it is one-to-one throughout the set $\delta B \times (1-\delta,1+\delta)$ for all small 
$\delta > 0$.  The image of this set under the map $Q$ is the given set, so the result now follows.
\finpf

We deduce the main result of this section.

\begin{thm}[Partial smoothness at disk matrices] \label{main-general}
Consider any matrix $X$ in $\Mn$.  Suppose there exists an invertible matrix $G \in \Mn$ such that
the largest eigenvalue of the matrix $\Phi(w^* X)$ is simple and independent of the scalar $w \in \T$, with a corresponding eigenvector $Gf(w)$.
Then around $X$, the set $\cA$ of disk matrices is an analytic manifold of codimension $2n$.  Its tangent and normal spaces at $X$ are
\begin{eqnarray*}
T_{\cA}(X) &=& \R X + \Big\{ Y \in \Mn : \ip{GF(w)G^*}{Y} = 0~ \mbox{for all}~ w \in \T \Big\} \\
N_{\cA}(X) &=& X^{\perp} \cap G \Big(\mbox{\rm span}\big\{ F(w) : w \in \T \big\}\Big) G^*.
\end{eqnarray*}
On $\cA$, around $X$, the numerical radius is analytic and its subdifferential is continuous and satisfies
\begin{eqnarray*}
\partial r (X) 
&=&  
G \Big(\mbox{\rm conv}\big\{ \|Gf(w)\|^{-2} F(w) : w \in \T \big\}\Big) G^* \\
\mbox{\rm par}\big(\partial r(X)\big)
&=&
N_{\cA}(X).
\end{eqnarray*}
Hence the numerical radius is partly smooth at the matrix $X$ relative to $\cA$.
\end{thm}

\pf
The case $n=1$ is an easy exercise, so we assume $n>1$.
Suppose first that the largest eigenvalue of $\Phi(w^* X)$ is $1$ for all $w$:  in other words,
the matrix $X$ has numerical radius $1$.
Using Theorem~\ref{normal} and the notation of its proof, to apply Proposition \ref{tool} we just need to check $X \not\in \T_{\cB}(X)$.  This follows from Lemma \ref{gradient}, which shows in fact 
$\RE\big( w^*g(w)^*Xg(w) \big) = 1$ for all $w \in \T$.  The properties of the subdifferential of the numerical radius follow from Theorem \ref{subdifferential-continuity}.  The result follows.

To prove the general case, define the positive scalar $\alpha = \frac{1}{r(X)}$, and apply the previous argument to the matrix $\alpha X$.  Since the set $\cA$ is closed under scalar multiplication, it is easy to check $N_{\cA}(\alpha X) = N_{\cA}(X)$.  Furthermore, since the numerical radius is positively homogeneous, it is also easy to check $\partial r(\alpha X) = \partial r (X)$.  This completes the proof.\finpf

\section{The Crouzeix parametrization} \label{parametrization}
According to Proposition \ref{eigenvalue}, a matrix $X \in \Mn$ belongs to the set $\cB$ of matrices with  field of values the unit disk if and only if the Hermitian matrix $\Phi(w^*X)$ has largest eigenvalue equal to one for all scalars $w \in \T$.  To guarantee  that $\cB$ is a manifold around $X$, Theorem \ref{normal} first made a modest refinement of the basic assumption:  that the largest eigenvalue is also simple for all $w \in \T$.  Secondly, however, we made a more opaque assumption:  that we can associate with this eigenvalue a corresponding eigenvector whose components are linearly independent polynomials in $w$ of degree less than $n$.  We now explore this assumption, using a result from \cite{Cro16}.  A somewhat related characterization appears in \cite[Thm 3.1]{MenOve05},

\begin{thm}[Crouzeix's disk matrix parametrization] \label{crouzeix}
A matrix in $\Mn$ has field of values the unit disk if and only if it has the form
$X = 2 (\sin B) U (\cos B)$ for some matrix $B \in \Hn$ satisfying 
\mbox{$\frac{\pi}{2}I \succcurlyeq B \succcurlyeq 0$} and some unitary $U \in \Mn$ such that the matrix 
$U(\cos B) - w(\sin B)$ is singular for all scalars  $w \in \T$.
\end{thm}

In this theorem statement, we can write 
\[
\sin B = V^* (\mbox{Diag}\, s) V ~~\mbox{and }~~ \cos B = V^* (\mbox{Diag}\, c) V, 
\]
for some unitary $V \in \Mn$ and nonnegative vectors $s,c \in \Rn$ satisfying $s_j^2 + c_j^2 = 1$ for $j=1,2,3,\ldots,n$.  With this notation, we have
\[
VXV^* ~=~ 2 (\mbox{Diag}\, s) (VUV^*) (\mbox{Diag}\, c)
\]
and the matrix $(VUV^*)(\mbox{Diag}\, c) - w (\mbox{Diag}\, s) )$ is singular for all $w \in \T$.  Hence we lose no essential generality in assuming, when convenient, that $V=I$, or in other words that $B$ is diagonal.

\begin{prop} \label{null}
For matrices $X$ as described in Theorem \ref{crouzeix}, the eigenspace corresponding to the largest eigenvalue of the matrix $\Phi(w^*X)$ is just the nullspace of the matrix $U(\cos B) - w(\sin B)$.
\end{prop}

\pf
By Theorem \ref{crouzeix} and Proposition \ref{eigenvalue}, we know 
$\lambda_{\max}\big(\Phi(w^*X)\big) = 1$.
Consider any vector $v$ in the nullspace of the matrix $U(\cos B) - w(\sin B)$.  We deduce
\begin{eqnarray*}
w^*(\sin B) U (\cos B) v & = &  (\sin B)^2 v  \\
w(\cos B) U^* (\sin B) v & = &  (\cos B)^2 v
\end{eqnarray*}
Adding shows $\Phi(w^*X)v = v$, as required.

Conversely, if $\Phi(w^*X)v = v$, then we have
\begin{eqnarray*}
\lefteqn{\big\| \big( U(\cos B)v - w(\sin B) v \big\|^2 ~=~
\|U(\cos B)v\|^2 + \|w(\sin B) v\|^2 } \\
& & \mbox{}
- \big(U(\cos B)v\big)^* \big(w(\sin B) v\big)
- \big(w(\sin B) v\big)^* \big(U(\cos B)v\big)  \\
& & = \|(\cos B)v\|^2 + \|(\sin B) v\|^2 - v^*\Phi(w^*X)v ~=~ 0, 
\end{eqnarray*}
as required.
\finpf

Any matrix with field of values the unit disk can be written in the form described in Theorem \ref{crouzeix}.  In that result, the matrix $U(\cos B) - w(\sin B)$ must be singular for all scalars 
$w \in \T$, or in matrix-analytic terminology, the matrix pencil \mbox{$U(\cos B) - w(\sin B)$} is singular.  Theorem \ref{normal} needs more:  in particular, the largest eigenvalue of the matrix 
$\Phi(w^*X)$ should be simple, or equivalently (by Proposition \ref{null}) the nullspace of $U(\cos B) - w(\sin B)$ (in the complex vector space $\Cn$) should be one-dimensional for all $w \in \T$.  Assuming this, the nullspace is in general spanned by some vector having all components polynomial in $w$, as can be seen quickly by reducing the matrix polynomial $U(\cos B) - w(\sin B)$ to Smith normal formal \cite{Gantmacher}.  However, Theorem \ref{normal} needs more precise knowledge of this eigenvector of $\Phi(w^*X)$.  Rather than delving into the theory of singular pencil canonical forms, we simply assume that the $n$ components of this vector can be taken to be linearly independent complex polynomials of degree less than $n$, or equivalently, that the vector can be written in the form $Gf(w)$, where the matrix $G \in \Mn$ is invertible.  This assumption will suffice for the examples we consider.

\begin{thm}[Partial smoothness and the Crouzeix parametrization] \label{main} \hfill \mbox{}
Suppose for some matrix $B \in \Hn$ satisfying $\frac{\pi}{2}I \succcurlyeq B \succcurlyeq 0$, and some unitary 
$U \in \Mn$, there exists an invertible $G \in \Mn$ such that
\[
\mbox{\rm Null} \big(U(\cos B) - w(\sin B)\big) ~=~ \C Gf(w) ~~ \mbox{for all scalars  $w \in \T$}.
\]
Then around the matrix $X = 2 (\sin B) U (\cos B)$ the set $\cA$ of disk matrices is an analytic manifold of codimension $2n$.  Its tangent and normal spaces at $X$ are
\begin{eqnarray*}
T_{\cA}(X) &=& \R X + \Big\{ Y \in \Mn : \ip{GF(w)G^*}{Y} = 0~ \mbox{for all}~ w \in \T \Big\} \\
N_{\cA}(X) &=& X^{\perp} \cap G \Big(\mbox{\rm span}\big\{ F(w) : w \in \T \big\}\Big) G^*.
\end{eqnarray*}
On $\cA$, around $X$, the numerical radius is analytic and its subdifferential is continuous and satisfies
\begin{eqnarray*}
\partial r (X) 
&=&  
G \Big(\mbox{\rm conv}\big\{ \|Gf(w)\|^{-2} F(w) : w \in \T \big\}\Big) G^* \\
\mbox{\rm par}\big(\partial r(X)\big)
&=&
N_{\cA}(X).
\end{eqnarray*}
In particular, the numerical radius is partly smooth at the matrix $X$ relative to the manifold $\cA$.
\end{thm}

\pf
This follows directly from Proposition \ref{null} and Theorem \ref{main-general}.
\finpf

We consider two special cases:  two-by-two matrices and superdiagonal matrices.

\subsection{The two-by-two case} \label{2by2}
As shown in \cite{Cro16}, a two-by-two matrix is a disk matrix if and only if it is unitarily similar to a multiple of the two-by-two Jordan block.  In Theorem \ref{main}, we could take
\[
U = 
\left[
\begin{array}{cc}
0 & 1 \\
1 & 0
\end{array}
\right]
~~\mbox{and}~~
B = 
\left[
\begin{array}{cc}
\frac{\pi}{2} & 0 \\
0 & 0
\end{array}
\right],
\]
so 
\[
X ~=~ 2 (\sin B) U (\cos B)
 ~=~ 2 
\left[
\begin{array}{cc}
1 & 0 \\
0 & 0
\end{array}
\right]
\left[
\begin{array}{cc}
0 & 1 \\
1 & 0
\end{array}
\right]
\left[
\begin{array}{cc}
0 & 0 \\
0 & 1
\end{array}
\right]
~=~
\left[
\begin{array}{cc}
0 & 2 \\
0 & 0
\end{array}
\right]
\]
and
\[
U(\cos B) - w(\sin B)
~=~
\left[
\begin{array}{cc}
0 & 1 \\
1 & 0
\end{array}
\right]
\left[
\begin{array}{cc}
0 & 0 \\
0 & 1
\end{array}
\right]
-
w
\left[
\begin{array}{cc}
1 & 0 \\
0 & 0
\end{array}
\right]
~=~
\left[
\begin{array}{cc}
-w & 1 \\
0 & 0
\end{array}
\right].
\]
Thus $G=I$ satisfies the assumptions of the theorem, so the numerical radius is partly smooth at $X$ relative to $\cA$.

\subsection{The superdiagonal case}
More generally, writing matrices in $\Mn$ (for $n > 2$) by partitioning both rows and columns into blocks of size 1, $n-2$, and 1, consider the case
\[
U = 
\left[
\begin{array}{ccc}
0 & e_1^* & 0 \\
0 & J & e_{n-2} \\
1 & 0 & 0
\end{array}
\right]
~~\mbox{and}~~
B = 
\left[
\begin{array}{ccc}
\frac{\pi}{2} & 0 & 0 \\
0 & \mbox{Diag}\,\theta & 0 \\
0 & 0 & 0
\end{array}
\right],
\]
where $J$ denotes the $(n-2)$-by-$(n-2)$ Jordan block, $e_1,e_{n-2} \in \R^{n-2}$ are the canonical unit vectors,
\[
J = 
\left[
\begin{array}{ccccc}
0 & 1 & 0 &  \ldots & 0 \\
0 & 0 & 1 &  \ldots & 0 \\
\vdots & \vdots &  \ddots & \ddots & \vdots \\
0 & 0 & 0 & \ddots & 1  \\
0 & 0 & 0 & \ldots & 0
\end{array}
\right],~~
e_1 = 
\left[
\begin{array}{c}
1 \\ 0 \\ \vdots \\ 0 \\ 0
\end{array}
\right],~~
e_{n-2} = 
\left[
\begin{array}{c}
0 \\ 0 \\ \vdots \\ 0 \\ 1
\end{array}
\right]
\]
 and the vector $\theta \in \R^{n-2}$ has entries in the interval $(0,\frac{\pi}{2})$.  Notice
\[
\sin B = 
\left[
\begin{array}{ccc}
1 & 0 & 0 \\
0 & \mbox{Diag}\,s & 0 \\
0 & 0 & 0
\end{array}
\right]
~~\mbox{and}~~
\cos B = 
\left[
\begin{array}{ccc}
0 & 0 & 0 \\
0 & \mbox{Diag}\,c & 0 \\
0 & 0 & 1
\end{array}
\right],
\]
where the entries $s_j,c_j \in (0,1)$ satisfy $s_j^2 + c_j^2 = 1$ for all $j=1,2,3,\ldots,n-2$, so
\begin{eqnarray*}
U(\cos B) - w(\sin B)
&=&
\left[
\begin{array}{ccc}
-w & c_1 e_1^* & 0 \\
0 & J (\mbox{Diag}\,c) - w(\mbox{Diag}\,s) & e_{n-2} \\
0 & 0 & 0
\end{array}
\right].
\end{eqnarray*}
A calculation shows that the matrix $G = \mbox{Diag}\,g$, where
\[
g_1 = 1 ~~\mbox{and}~~g_{j+1} = \frac{s_{j-1}}{c_j} g_j ~~(j=1,2,3,\ldots,n-1)
\]
(setting $s_0 = 1 = c_{n-1}$) satisfies the assumptions of Theorem \ref{main}.  Hence the numerical radius is partly smooth, relative to the manifold $\cA$, around the matrix 
\begin{eqnarray*}
X 
&=& 
2 (\sin B) U (\cos B) 
~=~ 2 
\left[
\begin{array}{ccc}
0 & c_1 e_1^* & 0 \\
0 & (\mbox{Diag}\,s) J (\mbox{Diag}\,c) & s_{n-2} e_{n-2} \\
0 & 0 & 0
\end{array}
\right],
\end{eqnarray*}
or more explicitly,
\bmye \label{super}
X ~=~
2\left[
\begin{array}{ccccc}
0 & s_0 c_1 & 0 &  \ldots & 0 \\
0 & 0 & s_1 c_2 &  \ldots & 0 \\
\vdots & \vdots &  \ddots & \ddots & \vdots \\
0 & 0 & 0 & \ddots & s_{n-2} c_{n-1}  \\
0 & 0 & 0 & \ldots & 0
\end{array}
\right].
\emye

For example, setting $s_j = \frac{1}{\sqrt{2}} = c_j$ for all $j=1,2,3,\ldots,n-2$ gives
\[
X =
\left[
\begin{array}{ccccccc}
0	& \sqrt{2} \\
	& \cdot	& 1 \\
	&		& \cdot	& \cdot \\
	&		&		& \cdot	& \cdot \\
	&		&		&		& \cdot	& 1 \\
	&		&		&		&		& \cdot	& \sqrt{2} \\
	&		&		&		&		&		& 0
\end{array}
\right],
\]
the Crabb matrix.  In fact, we have the following more general result.  We begin with a simple tool concerning real sequences.  For notational simplicity, and with no loss of generality, we present it for infinite sequences.

\begin{lem} \label{sequence}
For any strictly positive real sequence $b_1,b_2,b_3,\ldots$, there exists a decreasing sequence 
\[
+\infty = t_1 > t_2 > t_3 > \cdots > 0
\]
and a sequence of continuous, strictly increasing functions $h_j \colon [0,t_j) \to \R_+$, for all $j=1,2,3,\ldots$, satisfying 
\[
h_j(0) = 0,~~ h_j(t_{j+1}) = 1,~~ \lim_{t \uparrow t_j}h_j(t) = +\infty,
\]
$h_1(t) = b_1 t$ for all $t \in \R_+$, and 
\bmye \label{next}
h_{j+1}(t) = \frac{b_{j+1} t}{1 - h_j(t)} ~~ \mbox{for all}~ t \in [0,t_{j+1}).
\emye
\end{lem}

\pf
We construct the numbers $t_j$ and the functions $h_j$ recursively.  Given the continuous, strictly increasing function $h_j \colon [0,t_j) \to \R_+$ satisfying $h_j(0) = 0$ and 
$\lim_{t \uparrow t_j}h_j(t) = +\infty$ (as holds in particular when $j=1$), by the intermediate value theorem there exists a number $t_{j+1} \in (0,t_j)$ such that $h_j(t_{j+1}) = 1$.  Now the function 
$h_{j+1}$ defined by equation (\ref{next}) is continuous, strictly increasing and satisfies 
$h_{j+1}(0) = 0$ and $\lim_{t \uparrow t_{j+1}}h_{j+1}(t) = +\infty$, as required.
\finpf

We can now prove the following result.
\newpage

\begin{thm}[Partial smoothness at superdiagonal matrices] \hfill \mbox{} \\
At any matrix in $\Mn$ (with $n>1$) of the form
\bmye \label{general-super}
\left[
\begin{array}{ccccccc}
0	& a_1 \\
	& \cdot	& a_2 \\
	&		& \cdot	& \cdot \\
	&		&		& \cdot	& \cdot \\
	&		&		&		& \cdot	& a_{n-2}\\
	&		&		&		&		& \cdot	& a_{n-1} \\
	&		&		&		&		&		& 0
\end{array}
\right],
~~~\mbox{for nonzero $a_1,a_2,\ldots,a_{n-1} \in \C$},
\emye
the numerical radius is partly smooth relative to the set $\cA$ of disk matrices.
\end{thm}

\pf
We covered the case $n=2$ in Section \ref{2by2}, so assume $n>2$.
We claim that any matrix of the form (\ref{general-super}) 
is unitarily similar to a nonzero multiple of a matrix of the form (\ref{super}).  To see this, first note that, by diagonal unitary scaling we can first assume that each scalar $a_j$ is real and strictly positive.  We now seek numbers $s_j,c_j \in (0,1)$ satisfying $s_j^2 + c_j^2 = 1$ for all $j=1,2,3,\ldots,n-2$ and such that the vector
\[
[s_0c_1 ~~ s_1c_2 ~~ s_2c_3 ~\cdots~ s_{n-3}c_{n-2} ~~ s_{n-2}c_{n-1}]
\]
is a multiple of the vector
\[
[a_1 ~~ a_2 ~~ a_3 ~\cdots~ a_{n-2} ~~ a_{n-1}].
\]
As in the discussion before equation (\ref{super}), we set $s_0 = 1 = c_{n-1}$.
Squaring componentwise, we can write this problem as solving, for $c_j \in (0,1)$ and $t>0$, the equations
\[
c_1^2 = tb_1,~ (1-c_j^2)c_{j+1}^2 = tb_{j+1}~~(j=2,\ldots,n-2),~c_{n-1}=1,
\]
where $b_j = a_j^2$ for all $j$.  Using the notation of the lemma, these equations are equivalent to $c_j^2 = h_j(t)$ for $j \le n-1$, so are solvable:  the unique solution is $t=t_n$ and 
$c_j = \sqrt{h_j(t_n)}$ for $j \le n-1$.  

As we claimed, therefore, our matrix is unitarily similar to a nonzero multiple of a matrix of the form (\ref{super}).  As in the discussion of that example, the result now follows from Theorem \ref{main}.
\finpf

\noindent
We could, if desired, compute the numbers $c_j$ recursively as described in the proof of Lemma \ref{sequence}.

\section{The three-by-three case} \label{three-by}
We end with a thorough study of $3$-by-$3$ disk matrices.  In particular, we show that Theorem \ref{main-general} has applications beyond matrices with the sparsity pattern of a Jordan block and their unitary similarity transformations. 

Any disk matrix has a zero eigenvalue with algebraic multiplicity at least two, as shown in \cite[p.\ 438]{Cro16}.  Via Schur factorization, appropriately ordering the eigenvalues, any disk matrix in $\M^3$ is therefore unitarily similar to a matrix of the form
\bmye \label{ur}
\left[
\begin{array}{cc}
0 & S \\
0 & 0 
\end{array}
\right]
\quad
\mbox{with}
\quad
S = 2
\left[
\begin{array}{cc}
a & b \\
c & d 
\end{array}
\right]
\emye
for numbers $a,b,c,d \in \C$.  (The factor 2 is just for convenience, to simplify our later calculations.)  Furthermore, as we show below (in Theorem \ref{3disk}), the condition
\bmye \label{abcd}
c(|a|^2 + |d|^2) = -adb^*
\emye
must hold.  In other words, if we define two nonconvex cones of matrices,
\[
\Xi ~=~ \Big\{ \left[
\begin{array}{cc}
a & b \\
c & d 
\end{array}
\right]
 : 
 \mbox{condition (\ref{abcd}) holds}
\Big\}
\]
and
\[
\cE
~=~
\big\{
U [\begin{smallmatrix} 0 & S \\ 0 & 0 \end{smallmatrix}] U^* : 
U \in \M^3~ \mbox{unitary},~ S \in \Xi
\big\},
\]
then, in the space $\M^3$, the set of disk matrices $\cA$ is contained in $\cE$.
However, as we shall also see, this containment is strict.  

To organize the discussion, assuming condition (\ref{abcd}), exactly one of the following additional labeled conditions must hold:
\begin{eqnarray*}
a \ne 0 ~~\mbox{and}~~ d \ne 0 & & (\ne \ne) \\
a = 0 ~~\mbox{and}~~ d \ne 0 & & (= \ne) \\
a \ne 0 ~~\mbox{and}~~ d = 0 & & (\ne =) \\
a = 0 ~~\mbox{and}~~ d = 0  ~~\mbox{and}~~ 2|c| < |b| & & (= = <) \\
a = 0 ~~\mbox{and}~~ d = 0  ~~\mbox{and}~~ 2|c| = |b| & & (= = =) \\
a = 0 ~~\mbox{and}~~ d = 0  ~~\mbox{and}~~ 2|c| > |b| & & (= = >)
\end{eqnarray*}
Notice that the second and third conditions both entail $c=0$, by condition (\ref{abcd}), resulting in matrices of the form (\ref{ur}) that are nilpotent.

Based on these mutually exclusive conditions, we can partition the set $\Xi$ into corresponding disjoint subsets,
\bmye \label{partition}
\Xi ~~=~~ \Xi_{\ne \ne} ~\cup~ \Xi_{= \ne} ~\cup~ \Xi_{\ne =} ~\cup~ \Xi_{==<} ~\cup~ \Xi_{===} ~\cup~ \Xi_{==>},
\emye
and we can then write the set $\cE$ as a union of corresponding subsets,
\[
\cE ~~=~~ \cE_{\ne \ne} ~\cup~ \cE_{= \ne} ~\cup~ \cE_{\ne =} ~\cup~ \cE_{==<} ~\cup~ \cE_{===} ~\cup~ \cE_{==>}.
\]
This latter union is not, however, disjoint, as the following observation makes clear.

\begin{prop}
The sets $\cE_{= \ne}$ and $\cE_{\ne =}$ are identical and consist of the nilpotent matrices in the set 
$\cE_{==<}$.
\end{prop}

\pf
We simply observe, for example, that the matrices
\[
\left[
\begin{array}{ccc}
0 & 0 & b \\
0 & 0 & d \\
0 & 0 & 0
\end{array}
\right]
\quad
\mbox{and}
\quad
\left[
\begin{array}{ccc}
0 & 0 & |b|^2 + |d|^2 \\
0 & 0 & 0 \\
0 & 0 & 0
\end{array}
\right]
\]
are unitarily similar.  To see this, consider the unitary similarity transformation corresponding to  
$U = [\begin{smallmatrix} V & 0 \\ 0 & 1 \end{smallmatrix}]$, where the 2-by-2 unitary matrix $V$ maps 
the vector $[\begin{smallmatrix} b \\ d \end{smallmatrix}]$ to the vector 
$[\begin{smallmatrix} |b|^2 + |d|^2 \\ 0 \end{smallmatrix}]$. 
\finpf

With this result in mind, we can write
\[
\cE ~~=~~ \cE_{\ne \ne} ~\cup~ \cE_{==<} ~\cup~ \cE_{===} ~\cup~ \cE_{==>},
\]
and, as we shall see, this is indeed a partition into disjoint sets, and one well suited to understanding the set $\cA$ of disk matrices and the numerical radius $r$.  Specifically, we devote this section to showing the following result.  Recall that $\Omega$ is the open subset of $\M^3$ defined in equation~(\ref{Omega}).

\begin{thm} \label{summary}
In the 18-dimensional space $\M^3$, the set $\cE_{\ne \ne}$ is a 12-dimensional manifold.  Its closure is the set of disk matrices, which can be partitioned into three disjoint sets:
\[
\cA ~=~ \mbox{\rm cl}\, \cE_{\ne \ne} ~=~
\cE_{\ne \ne} ~\cup~ \cE_{==<} ~\cup~ \cE_{===}.
\]
These sets satisfy
\[
\cA \cap \Omega ~=~ \cE_{\ne \ne} ~\cup~ \cE_{==<}.
\]
The subdifferential of the numerical radius $\partial r(\cdot)$ has dimension 6 throughout 
$\cE_{\ne \ne}$ and 4 throughout $\cE_{==<}$.  The numerical radius is partly smooth at a disk matrix $X$ relative to $\cA$ if and only if $X \in \cE_{\ne \ne}$.
\end{thm}

We begin the proof with some simple algebraic observations about the partition~(\ref{partition}).

\begin{prop} \label{closure}
The following relationships hold:
\[
\mbox{\rm cl}\, \Xi_{\ne \ne} = \Xi \setminus \Xi_{==>}
\quad \mbox{and} \quad
\mbox{\rm cl}\, \cE_{\ne \ne} = \cE \setminus \cE_{==>}.
\]
\end{prop}

\pf
For any number $b \in \C$ and nonzero $d \in \C$, we have
\[
\lim_{\epsilon \to 0}
\left[
\begin{array}{cc}
\epsilon & b \\
-\frac{db^*\epsilon}{\epsilon^2 + |d|^2} & d
\end{array}
\right]
~=~
\left[
\begin{array}{cc}
0 & b \\
0 & d
\end{array}
\right],
\]
so
\bmye \label{=ne}
\Xi_{=\ne} ~\subset~ \mbox{cl}(\Xi_{\ne\ne}),
\emye
and similarly $\Xi_{\ne =} \subset \mbox{cl}(\Xi_{\ne\ne})$.
On the other hand, for any numbers $c,d \in \C$ satisfying $2|c| < |b|$, the real number
\[
t = \frac{1}{2} \Big( \sqrt{1 - \Big( \frac{2|c|}{|b|} \Big)^2} - 1 \Big)
\]
satisfies
\[
t^2 + t  + \frac{|c|^2}{|b|^2} = 0,
\]
so for $\epsilon > 0$ we have
\[
\left[
\begin{array}{cc}
\epsilon c & b \\
c & \epsilon b t
\end{array}
\right]
~\in~ 
\left\{
\begin{array}{ll} 
\Xi_{\ne \ne} & (c \ne 0) \\
\Xi_{= \ne} & (c = 0).
\end{array}
\right.
\]
This matrix converges to the matrix $[\begin{smallmatrix} 0 & b \\ c & 0 \end{smallmatrix}]$ as 
$\epsilon \downarrow 0$, so using property (\ref{=ne}) we deduce
$\Xi_{==<} \subset \mbox{cl}(\Xi_{\ne\ne})$.
Since we obviously have
\bmye \label{obvious}
\Xi_{===} ~\subset~ \mbox{cl}(\Xi_{==<}),
\emye
we deduce $\mbox{\rm cl}( \Xi_{\ne \ne} ) \supset \Xi \setminus \Xi_{==>}$.
The opposite inclusion is easy: when $a$ and $d$ are nonzero, condition (\ref{abcd}) implies $2|c| \le |b|$, since 
\[
\frac{|a|^2 + |d|^2}{|ad|} \ge 2.
\]
The second equation is a simple consequence, using the compactness of the set of unitary matrices.
\finpf

We next introduce a simple tool.

\begin{lem} \label{previous}
Consider a matrix $X$ unitarily similar to a matrix of the form (\ref{ur}). The characteristic polynomial of the matrix $\Phi(w^*X)$, namely $\det\big(\Phi(w^* X) - \lambda I\big)$, is
\[
-\lambda^3 + 2\lambda^2 \ip{w}{c}  + (|a|^2 + |d|^2 + |b|^2) \lambda 
+ 2\ip{w}{adb^* - c|b|^2}.
\]
\end{lem}

\pf
If $X$ is unitarily similar to a matrix $E$ in the form (\ref{ur}), then the characteristic polynomial of $\Phi(w^*X)$ is identical with that of $\Phi(w^*E)$, namely
\[
\det\big(\Phi(w^* E) - \lambda I\big)
~=~
\det\left[
\begin{array}{ccc} \label{upper-right}
-\lambda 	& w^*a 					& w^*b \\
wa^* 		& wc^* + w^*c - \lambda & w^*d \\
wb^* 		& wd^* 					& -\lambda
\end{array}
\right],
\]
and the result follows.
\finpf

We can now characterize disk matrices in $\M^3$.

\begin{thm} \label{3disk}
The set of disk matrices in $M^3$ is $\cE \setminus \cE_{==>}$.
\end{thm}

\pf
We need to show that a matrix $X \in \M^3$ is a disk matrix if and only if it is unitarily similar to a matrix of the form $[\begin{smallmatrix} 0 & S \\ 0 & 0 \end{smallmatrix}]$ for a matrix $S \in \Xi \setminus \Xi_{==>}$.  To this end, consider any $X$ with field of values a disk of radius $r \ge 0$.  As we noted, $X$ is  unitarily similar to a matrix of the form (\ref{ur}), so by Lemma \ref{previous} the characteristic polynomial of $\Phi(w^*X)$ is  
\[
-\lambda^3 + 2\lambda^2 \ip{w}{c}  + (|a|^2 + |d|^2 + |b|^2) \lambda 
+ 2\ip{w}{adb^* - c|b|^2}.
\]
By Proposition \ref{eigenvalue}, for all $w \in \T$ we have
\[
-r^3 + 2r^2 \ip{w}{c}  + (|a|^2 + |d|^2 + |b|^2) r 
+ 2\ip{w}{adb^* - c|b|^2} ~=~ 0.
\]
By considering $w = \pm 1$ and adding we deduce
\[
r = \sqrt{|a|^2 + |d|^2 + |b|^2}.
\]
Substituting into the previous equation gives
\[
\ip{w}{ c(|a|^2 + |d|^2) + adb^*}
 ~=~ 0
\]
for all $w \in \T$, and equation (\ref{abcd}) follows.  

If in addition we have $S \in \Xi_{==>}$, then 
$r = |b|$.  Our previous calculation shows
\begin{eqnarray*}
\lefteqn{
\det\big(\Phi(w^* E) - \lambda I\big)
} \\
&=& 
-\lambda^3 + 2\lambda^2 \ip{w}{c}  + |b|^2 \lambda - 2|b|^2\ip{w}{c} \\
&=&
-(\lambda - 2\ip{w}{c}) \big(\lambda^2 - |b|^2) \big),
\end{eqnarray*}
so the matrix $\Phi(w^* X)$ has eigenvalues $2\ip{w}{c}$ and $\pm |b|$.  Setting $w = \frac{1}{|c|}c$ now gives $r=\max(2|c|,|b|) = 2|c|$, a contradiction.

Conversely, consider any matrix $E$ of the form (\ref{ur}), where $S \in \Xi$ but $S \notin \Xi_{==>}$.  We want to prove 
\[
\lambda_{\max}\big(\Phi(w^*E)\big) = \sqrt{|a|^2 + |d|^2 + |b|^2} ~~\mbox{for all}~ w \in \T,
\]
which will imply that 
$E$, along with any matrix that is unitarily similar to $E$, is a disk matrix.  Lemma \ref{previous} shows
\begin{eqnarray*}
\lefteqn{
\det\big(\Phi(w^* E) - \lambda I\big)
} \\
&=& 
-\lambda^3 + 2\lambda^2 \ip{w}{c}  + (|a|^2 + |d|^2 + |b|^2) \lambda 
- 2(|a|^2 + |d|^2 + |b|^2)\ip{w}{c} \\
&=&
-(\lambda - 2\ip{w}{c}) \big(\lambda^2 - (|a|^2 + |d|^2 + |b|^2) \big),
\end{eqnarray*}
so the matrix $\Phi(w^* E)$ has eigenvalues $2\ip{w}{c}$ and $\pm \sqrt{|a|^2 + |d|^2 + |b|^2}$.  The result now follows immediately in the case when $a=0=d$, since then by assumption $2|c| \le |b|$.
It also follows if either $a$ or $d$ is nonzero, since equation (\ref{abcd}) then implies
\[
2|c| < \sqrt{|a|^2 + |d|^2 + |b|^2}
\]
To see this, note the equivalent form
\[
2|c|(|a|^2 + |d|^2) < (|a|^2 + |d|^2)\sqrt{|a|^2 + |d|^2 + |b|^2}
\]
or equivalently
\[
2|a||d||b| < (|a|^2 + |d|^2)\sqrt{|a|^2 + |d|^2 + |b|^2}.
\]
Squaring gives
\[
4|a|^2|d|^2|b|^2 < (|a|^2 + |d|^2)^2(|a|^2 + |d|^2 + |b|^2)
\]
or equivalently the valid inequality
\[
0 < (|a|^2 - |d|^2)^2 |b|^2 + (|a|^2 + |d|^2)^3.
\]
This completes the proof of the converse direction.
\finpf

Examining the proof, we have in particular shown the following result.  Notice that, having already discarded the subset $\Xi_{==>}$, this reveals the special role of the subset 
$\Xi_{===}$.

\begin{lem} \label{general-disk}
A matrix $E$ of the form (\ref{ur}) is a disk matrix if and only if the matrix $S$ lies in the set 
$\Xi \setminus \Xi_{==>}$.  In that case, $E$ has numerical radius
\[
\lambda_{\max}\big(\Phi(w^*E)\big) = \sqrt{|a|^2 + |d|^2 + |b|^2} ~~\mbox{for all}~ w \in \T.
\]
Furthermore this largest eigenvalue is simple for all $w \in \T$ if and only if $S \not\in \Xi_{===}$.
\end{lem}

Discarding the set $\Xi_{===}$ rules out the case of the zero matrix, so by scaling we can focus our attention on matrices with field of values the unit disk.  In this context, we next further discard the set $\Xi_{==<}$.

\begin{lem} \label{three-a}
A matrix $E$ of the form (\ref{ur}) has field of values the unit disk if and only if the matrix $S$ lies in the set $\Xi \setminus \Xi_{==>}$ and $|a|^2 + |d|^2 + |b|^2 = 1$.  Assume this is the case.  Then 
\[
\lambda_{\max}\big(\Phi(w^*E)\big) = 1 ~~\mbox{for all}~ w \in \T.
\]
Furthermore, this largest eigenvalue is simple for all $w \in \T$ if and only if we have
$S \not\in \Xi_{===}$.  If we further assume that $a$ and $d$ are not both zero 
(so $S \not\in \Xi_{==<}$), then
a corresponding unit eigenvector is
\[
 \frac{1}{\sqrt{2(1-|b|^2)(1 - 2\ip{w}{c})}}
\left[
\begin{array}{c}
w^*a + bd^* \\
1 - |b|^2 \\
b^*a + wd^*
\end{array}
\right].
\]
\end{lem}

\pf
By our previous calculations we know, for all $w \in \T$, that the matrix
\[
\Phi(w^*E) - I
~=~
\left[
\begin{array}{ccc} 
-1 	& w^*a 					& w^*b \\
wa^* 		& wc^* + w^*c - 1 & w^*d \\
wb^* 		& wd^* 					& -1
\end{array}
\right]
\]
has nullity one.  It is easy to verify
\[
\left[
\begin{array}{ccc} 
-1 	& w^*a 					& w^*b \\
wb^* 		& wd^* 					& -1
\end{array}
\right]
\left[
\begin{array}{c}
w^*a + bd^* \\
1 - |b|^2 \\
b^*a + wd^*
\end{array}
\right]
~=~
\left[
\begin{array}{c}
0 \\
0
\end{array}
\right],
\]
so the given vector is indeed an eigenvector, and the square of its norm is
\[
\frac{|a|^2 + |b|^2|d|^2 + 2\ip{w}{adb^*} + (1-|b|^2)^2 + |a|^2|b|^2 + |d|^2 + 2\ip{w}{adb^*}}
{2(1-|b|^2)(1 - 2\ip{w}{c})}
~=~ 1
\]
as required.
\finpf

We are now ready to consider partial smoothness.

\begin{thm} \label{three-final}
Consider any matrix $E$ of the form (\ref{ur}), for a matrix $S$ in the set  
$\Xi_{\ne \ne}$.  Then $E$ belongs to the set $\cA$ of disk matrices.  Around $E$, furthermore, 
$\cA$ is a 12-dimensional analytic manifold in the 18-dimensional space $\M^3$, with normal space at $E$ given by
\[
N_{\cA}(E) ~=~ E^{\perp} \cap G \Big(\mbox{\rm span}\big\{ F(w) : w \in \T \big\}\Big) G^*,
\]
where $G$ is the invertible matrix
\[
G = 
\left[
\begin{array}{ccc}
a 	& bd^* 		& 0 \\
0 	& 1-|b|^2 	& 0 \\
0	& ab^*		& d^*
\end{array}
\right].
\]
The numerical radius is partly smooth at $E$ relative to $\cA$, and its subdifferential there, $\partial r(E)$, has dimension 6.
\end{thm}

\pf
As in the proof of Theorem \ref{main-general}, we lose no generality in first scaling the matrix $E$ to ensure $|a|^2 + |d|^2 + |b|^2 = 1$.  Then, by Lemma \ref{three-a}, the eigenvalue
\[
\lambda_{\max}\big(\Phi(w^*E)\big) = 1 
\]
is simple for all scalars $w \in \T$, with a corresponding eigenvector 
\[
Gf(w) =
\left[
\begin{array}{c}
aw + bd^*w^2 \\
(1 - |b|^2)w^2 \\
b^*aw^2 + d^*w^3
\end{array}
\right].
\]
The result now follows by Theorem \ref{main-general}.
\finpf

\begin{cor}
The semi-algebraic set $\cA$ of all disk matrices in $\M^3$ has dimension~12.
\end{cor}

\pf
By Proposition \ref{closure} and Theorem \ref{3disk}, any disk matrix is a limit of matrices $X$ unitarily similar to those of the form considered in Theorem \ref{three-final}.  That result shows that the set of such matrices $X$ is a semi-algebraic manifold of dimension 12, so its closure also has dimension 12.
\finpf

By contrast, we have the following result.

\begin{prop} \label{failure}
Any matrix $X \in \cE_{= = <}$ is a disk matrix:  in fact $X \in \cA \cap \Omega$.  However,
the subdifferential of the numerical radius $\partial r(X)$ has dimension 4.  Consequently, the numerical radius is not partly smooth at $X$ relative to any subset of $\cA$.
\end{prop}

\pf
After a unitary similarity transformation, we can suppose
\[
X ~=~
2 \left[
\begin{array}{ccc}
0 & 0 & b \\
0 & c & 0 \\
0 & 0 & 0
\end{array}
\right],
\]
with $2|c| < |b|$.  The claim $X \in \cA \cap \Omega$ follows from Lemma \ref{three-a}.
By Lemma \ref{general-disk}, for all scalars 
$w \in \T$, the eigenvalue
\[
\lambda_{\max}\big(\Phi(w^* X) \big) = |b| 
\]
is simple, and a quick calculation shows that 
a corresponding unit eigenvector is
\[
 \frac{1}{\sqrt{2}}
 \big[1 ~~ 0 ~~ uw^*\big]^*
\]
where $u = \frac{b}{|b|}$.
Now Theorem \ref{subdifferential-continuity} implies that the numerical radius has subdifferential at $X$ given by
\[
\partial r(X) ~=~
\frac{1}{2} \mbox{conv} \left\{ 
\left[
\begin{array}{ccc}
w 		& 0 		& u 		\\
0	& 0	& 0 	\\
u^* w^2 & 0 	& w
\end{array}
\right]
: w \in \T 
\right\}.
\]
We can rewrite the matrix in this equation as the image of the vector 
$[1 ~ w ~ w^2]$ under a one-to-one linear map,
so by Corollary \ref{dimension} we see that $\partial r(X)$ has dimension~$4$.  Finally, by definition, any manifold with respect to which the numerical radius is partly smooth at $X$ must have dimension $18-4 = 14$, and hence cannot be a subset of the 12-dimensional set $\cA$.
\finpf

We can summarize some of our observations in the following refinement of Theorem \ref{3disk}, giving the final ingredient in the proof of Theorem \ref{summary}.

\begin{cor}
For any matrix $X$ in the set of disk matrices $\cA$ in $\M^3$, the numerical radius is partly smooth at $X$ relative to $\cA$ if and only if $X \in \cE_{\ne \ne}$.
\end{cor}

\pf
We can suppose that $X$ is unitarily similar to a matrix $E=[\begin{smallmatrix} 0 & 2S \\ 0 & 0 \end{smallmatrix}]$ for a matrix $S \in \Xi \setminus \Xi_{==>}$.  We claim that the numerical radius is  partly smooth at $E$ relative to $\cA$ if and only if $S \in \Xi_{\ne \ne}$.
In the light of Theorem \ref{three-final} and Proposition~\ref{failure},
the only case left to consider is $S \in \Xi_{===}$.  In the case $S \in \Xi_{\ne \ne}$ we know 
$\dim \partial r(E) = 6$, whereas in the case $S \in \Xi_{==<}$ we have $\dim \partial r(E) = 4$.  But when $S \in \Xi_{===}$, our matrix can be written both as a limit of matrices $E \in \cA$ with $S \in \Xi_{\ne \ne}$ or alternatively with $S \in \Xi_{==<}$, by Proposition \ref{closure} and equation (\ref{obvious}), so $\dim \partial r$ is not continuous there relative to $\cA$, contradicting partial smoothness.
\finpf

It is interesting to consider the impact of these results on optimization problems involving the numerical radius $r$ on $\M^3$.  For example, by Theorems \ref{identification} and \ref{three-final},  the prox operator $\mbox{prox}_r$ commonly maps a random matrix $Y$ to a disk matrix, which is 
therefore unitarily similar to a matrix $E=[\begin{smallmatrix} 0 & 2S \\ 0 & 0 \end{smallmatrix}]$ for a matrix $S \in \Xi \setminus \Xi_{==>}$.  In particular, the case $S \in \Xi_{\ne \ne}$
occurs with positive probability, and this property is stable under perturbations to $Y$.  

Worth noting is that the case $S \in \Xi_{===}$ can also occur stably with positive probability.  For example, the prox operator $\mbox{prox}_r$ maps any small matrix to the zero matrix. 

On the other hand, the case $S\in \Xi_{= = <}$ covered by Proposition \ref{failure} (as well as the cases $S \in \Xi_{= \ne} \cup \Xi_{\ne =}$) do not display this same stable behavior.  For almost all matrices $Y$, the numerical radius, being semi-algebraic, will indeed be partly smooth around $\mbox{prox}_r(Y)$ relative to some manifold $\cM$ (see \cite{gen}).  However, even if $\mbox{prox}_r(Y)$ were a disk matrix of one of these three forms, $\cM$ could not be a subset of $\cA$ locally, so certain small perturbations to $Y$ would produce matrices other than disk matrices.  

\subsubsection*{The Crouzeix ratio in the 3-by-3 case}
It is interesting to compare how two $3\times 3$ families given in \cite{Cro16} fit into the discussion of this section. Let 
\begin{equation*}
  Y(\xi)= \left [ \begin{array}{ccc} 
             0 & 0 & 2\\
             0 & \xi & 0 \\
             0 & 0 & 0
             \end{array} \right ], \ \xi\in\C, \ |\xi| \leq 1   
\end{equation*}
and
\begin{equation*}
  Z(\phi,\psi)=  e^{i \psi} \left [ \begin{array}{ccc} 
             0 & \sqrt{2} \cos\phi  & 2\ \sin\phi\\
             0 & -\sin \phi &  \sqrt{2}\cos\phi \\
             0 & 0 & 0
             \end{array} \right ], \ \phi\in\left[0, \frac{\pi}{2}\right], \ \psi\in\R. 
\end{equation*}
Crouzeix proved in \cite[Lemma 10.2]{Cro16} that these matrices, along with their unitary similarity
transformations, are the only $3 \times 3$ matrices $X$ with field of values the unit disk and
\[
     \inf\left \{\frac{\max\{|p(\lambda)|: |\lambda| \le 1\}}{\|p(X)\|_2}: p \text{~is a polynomial of any degree~} \right \}  = \frac{1}{2}.
\]
Note that the ratio in (\ref{cro-ratio}) is a special
case of the one in the equation above, with $p(\lambda)=\lambda^{n-1}$. Indeed, since permuting the
rows and columns of $Y(\xi)$ results in a block diagonal matrix with one block equal to the scaled $2\times 2$
Jordan block discussed in Section~\ref{examples}, and a second block equal to the scalar $\xi$ with $|\xi|\leq 1$,
when $X=Y(\xi)$ the infimum in the equation above is attained by $p(\lambda)=\lambda$. When $\phi=\frac{\pi}{2}$,
the matrix $Z(\phi,\psi)$ can also, by permuting rows and columns, be transformed to essentially the same block diagonal form,
while when $\phi=0$, $Z(\phi,\psi)$ is a unimodular scalar times the scaled $3\times 3$
Jordan block discussed in Section~\ref{examples}, and the infimum is attained by  $p(\lambda)=\lambda^2$. In general, however, 
for $X=Z(\phi,\psi)$ the infimum is not attained (although it is if $p$ is allowed to be any analytic function).

We can apply the results in this section directly to the matrices $Y(\xi)$ and $Z(\phi,\psi)$.  In particular, the numerical radius is partly smooth relative to the set of disk matrices at the matrix $Z(\phi,\psi)$ if and only if  $\phi\not=\frac{\pi}{2}$, and is not so at the matrix $Y(\xi)$.

\section{Concluding Remarks}
In this paper, we have analyzed the partly smooth geometry of the numerical radius
around disk matrices --- those whose fields of values are disks in the complex plane centered
at the origin. We have argued that such matrices arise naturally as solutions of optimization problems
involving the numerical radius, supported by both numerical experiments and subdifferential analysis.
Even in the case of matrices of order three, the details are surprisingly subtle.  Our development blends matrix and variational analysis, but makes no use of semidefinite programming characterizations of the numerical radius like Theorem \ref{mathias}:  the implication of this characterization for partial smoothness is an interesting open question.

\bibliography{adrian}
\bibliographystyle{plain}

\end{document}